# The Marshall-Olkin-Kumaraswamy-G family of distributions

Laba Handique and Subrata Chakraborty*
Department of Statistics, Dibrugarh University
Dibrugarh-786004, India
*Corresponding author: Email: subrata_stats@dibru.ac.in

(August 21, 2016)

**Abstract**

A new family of continuous distribution is proposed by using Kumaraswamy-*G* (Cordeiro and de Castro, 2011) distribution as the base line distribution in the Marshal-Olkin (Marshall and Olkin, 1997) construction. A number of known distributions are derived as particular cases. Various properties of the proposed family like formulation of the pdf as different mixture of exponentiated baseline distributions, order statistics, moments, moment generating function, Renyi entropy, quantile function and random sample generation have been investigated. Asymptotes, shapes and stochastic ordering are also investigated. The parameter estimation by methods of maximum likelihood, their large sample standard errors and confidence intervals and method of moment are also presented. Two members of the proposed family are compared with corresponding members of Kumaraswamy-Marshal-Olkin-*G* family (Alizadeh *et al.*, 2015) by fitting of two real life data sets.

---

*Preprint*

## 1. Introduction

One of the preferred area of research in the filed of the probability distribution is that of generating new distributions starting with a base line distribution by adding one or more additional parameters notable among them are Azzalini's skewed family (Azzalini, 1985), Marshall-Olkin extended ($MOE$) family (Marshall and Olkin, 1997), exponentiated family ($EF$) of distributions (Gupta *et al.*, 1998), or by composite methods of combining two or more known competing distribution through transformations like beta-generated ($beta-G$) family (Eugene *et al.*, 2002; Jones 2004), gamma-generated ($G-G$) families (Zografos Balakrishnan 2009, Ristic and Balakrishnan, 2012), Kumaraswamy - $G$ ($Kw-G$) family (Cordeiro and de Castro, 2011 ; Nadarajah *et al.*, 2012; Hussain, 2013), McDonald-$G$ ($Mc-G$) (Alexander *et al.*, 2012), beta extended-$G$ family (Cordeiro *et al.*, 2012), Kumaraswamy -beta generalized family (Pescim *et al.*, 2012), exponentiated transformed transformer family ($ET-X$), (Alzaghal *et al.* 2013), exponentiated generalized ($Exp-G$) family (Cordeiro *et al.*, 2013), geometric exponential-Poisson family (Nadarajah *et al.*, 2013a), truncated-exponential skew symmetric family (Nadarajah *et al.*, 2013b), logistic-generated ($Lo-G$) family (Torabi and Montazari, 2014) and Kumaraswamy Marshal - Olkin family (Alizadeh *et al.*, 2015). While the additional parameter(s) bring in more flexibility at the same time they also complicates the mathematical form of the resulting distribution, often considerably enough to render it not amenable to further analytical and numerical manipulations. But with the advent of sophisticated powerful mathematical and statistical softwares unlike in past now a days more and more complex distributions are getting accepted as viable models of data analysis. Tahir *et al.* (2015) provided a detail account of how new families of univariate continuous distributions are being generated through introduction of additional parameters.

### 1.1 Formulas and notations

Here we list some formulas to be used in the subsequent sections of this article.

If $T$ is a continuous random variable with pdf, $f(t)$ and cdf $F(t)=P[T\le t]$, then its

Survival function (sf): $\overline{F}(t)=P[T>t]=1-F(t)$,

Hazard rate function (hrf): $h(t)=f(t)/\overline{F}(t)$,

Reverse hazard rate function (rhrf): $r(t)=f(t)/F(t)$,

Cumulative hazard rate function (chrf): $H(t)=-\log[\overline{F}(t)]$.

$(p,q,r)^{\text{th}}$ Power Weighted Moment (PWM): $\Gamma_{p,q,r}=\int_{-\infty}^{\infty} t^p[F(t)]^q[1-F(t)]^r f(t)\,dt$,

Entropy: The Rényi entropy of $T$ is given by $I_R(\delta) = (1-\delta)^{-1} \log\left( \int_{-\infty}^{\infty} f(t)^{\delta} \, dt \right)$.

**1.2 Marshall-Olkin Extended $(MOE)$ family of distributions**

Starting with a given baseline distribution with probability density function (pdf) $f(t)$, cumulative distribution function (cdf) $F(t)$, Marshall and Olkin (1997) proposed a new flexible semi parametric family of distributions and defined a new survival function (sf) $\overline{F}^{MO}(t)$ by introducing an additional parameter $\alpha > 0$. The sf $\overline{F}^{MO}(t)$ of the MOE family of distributions is defined by

$$\overline{F}^{MO}(t) = \frac{\alpha \overline{F}(t)}{1-\overline{\alpha}\,\overline{F}(t)} = \frac{\alpha \overline{F}(t)}{F(t)+\alpha \overline{F}(t)} = \frac{\alpha\,[1-F(t)]}{\alpha+\overline{\alpha}\,F(t)} \tag{1}$$

where $-\infty < t < \infty, \alpha > 0$ and $\overline{\alpha} = 1-\alpha$. The parameter $'\alpha'$ is known as the tilt parameter as since the hrf of the new family is shifted below (above) for $\alpha \geq 1 (0 < \alpha \leq 1)$ the hrf of the base line distribution (Nanda and Das, 2012). That is for all $t \geq 0, h^{MO}(t) \leq h(t)$ when $\alpha \geq 1$, and $h^{MO}(t) \geq h(t)$ when $0 < \alpha \leq 1$, where $h^{MO}(t)$ and $h(t)$ are the hrf's of the $MOE$ and baseline distributions respectively. Now

$$F^{MO}(t) = 1-\overline{F}^{MO}(t) = 1-\frac{\alpha \overline{F}(t)}{1-\overline{\alpha}\,\overline{F}(t)} = \frac{1-(\alpha+\overline{\alpha})\overline{F}(t)}{1-\overline{\alpha}\,\overline{F}(t)} = \frac{F(t)}{1-\overline{\alpha}\,\overline{F}(t)} = \frac{F(t)}{1-(1-\alpha)\overline{F}(t)}$$

$$= \frac{F(t)}{F(t)+\alpha \overline{F}(t)} \text{ or } \frac{1-\overline{F}(t)}{\alpha+\overline{\alpha}\,F(t)} = \frac{F(t)}{\alpha+\overline{\alpha}\,F(t)} \tag{2}$$

and $$f^{MO}(t) = \frac{d}{dt} F^{MO}(t) = \frac{d}{dt}\left[\frac{F(t)}{\alpha+\overline{\alpha}\,F(t)}\right] = \frac{\alpha\, f(t)}{[\alpha+\overline{\alpha}\,F(t)]^2} \text{ or } \frac{\alpha\, f(t)}{[1-\overline{\alpha}\,\overline{F}(t)]^2} \tag{3}$$

Where $-\infty < t < \infty, \alpha > 0$ and $\overline{\alpha} = 1-\alpha$. If $\alpha = 1$, then we have $\overline{F}^{MO}(t) = \overline{F}(t)$. Other reliability measures like the hrf, rhrf and chrf associated with (1) are

$$h^{MO}(t) = \frac{f^{MO}(t)}{\overline{F}^{MO}(t)} = \frac{\alpha\, f(t)}{[1-\overline{\alpha}\,\overline{F}(t)]^2} \Bigg/ \frac{\alpha \overline{F}(t)}{1-\overline{\alpha}\,\overline{F}(t)} = \frac{f(t)}{\overline{F}(t)} \frac{1}{[1-\overline{\alpha}\,\overline{F}(t)]} = \frac{h(t)}{1-\overline{\alpha}\,\overline{F}(t)}$$

$$r^{MO}(t) = \frac{f^{MO}(t)}{F^{MO}(t)} = \alpha \frac{f(t)}{F(t)} \frac{1}{[1-\overline{\alpha}\,\overline{F}(t)]} = \frac{\alpha\, r(t)}{1-\overline{\alpha}\,\overline{F}(t)} \text{ or } \frac{\alpha\, r(t)}{\alpha+\overline{\alpha}\,F(t)} \quad \text{and}$$

$$H^{MO}(t) = -\log\left[\frac{\alpha \overline{F}(t)}{1-\overline{\alpha}\,\overline{F}(t)}\right] \text{ or } -\log\left\{\frac{\alpha[1-F(t)]}{\alpha+\overline{\alpha}\,F(t)}\right\} \quad \text{respectively.}$$

Where $h(t)$ and $r(t)$ is the hrf and rhrf of the baseline distribution.

It is obvious that many new families can be derived from Marshall-Olkin set up by considering different base line distribution in the equation (1). These new families are usually termed as Marshall-Olkin extended *F* distribution. For details see Tahir *et al.* (2015).

Jayakumar and Mathew (2008) generalized the Marshall-Olkin set up by exponentiating the Marshall-Olkin survival function

$$\bar{F}^{GMO}(t)=\left[\frac{\alpha\,\bar{F}(t)}{1-\bar{\alpha}\,\bar{F}(t)}\right]^{\theta}.$$

This method is called method of Lehmann alternative 1 due to Lehmann (1953)

Tahir *et al.* (2015) propose another generalization through Lehmann alternative 2 due to Lehmann (1953) by exponentiating the Marshall-Olkin cdf as

$$F^{GMO}(t)=\left[1-\frac{\alpha\,\bar{F}(t)}{1-\bar{\alpha}\,\bar{F}(t)}\right]^{\theta}.$$

[for more on Marshall - Olkin family check Barakat *et al.* (2009), Jose (2011), Krishna (2011), Barreto-Souza *et al.* (2013), Cordeiro et al. (2014a), Alizadeh *et al.*, (2015)]

### 1.3 Kumaraswamy-G $(Kw-G)$ family of distributions

For a baseline cdf $G(t)$ with pdf $g(t)$, Cordeiro and de Castro (2011) defined $Kw-G$ distribution with cdf and pdf

$$F^{KwG}(t)=1-[1-G(t)^{a}]^{b} \tag{4}$$

and

$$f^{KwG}(t)=a\,b\,g(t)G(t)^{a-1}[1-G(t)^{a}]^{b-1} \tag{5}$$

Where $t>0,\ g(t)=\frac{d}{dt}G(t)$ and $a>0, b>0$ are shape parameters in addition to those in the baseline distribution which partly govern skewness and variation in tail weights. For a lifetime random variable $'t'$, the sf, hrf, rhrf and chrf for distribution in (4) are given respectively by

$$\bar{F}^{KwG}(t)=1-F^{KwG}(t)=[1-G(t)^{a}]^{b}$$

$$h^{KwG}(t)=f^{KwG}(t)\big/\bar{F}^{KwG}(t)=a\,b\,g(t)G(t)^{a-1}[1-G(t)^{a}]^{b-1}\big/[1-G(t)^{a}]^{b}$$

$$=a\,b\,g(t)G(t)^{a-1}[1-G(t)^{a}]^{-1}$$

$$r^{KwG}(t)=f^{KwG}(t)\big/\bar{F}^{KwG}(t)=a\,b\,g(t)G(t)^{a-1}[1-G(t)^{a}]^{b-1}\big/1-[1-G(t)^{a}]^{b}$$

$$=a\,b\,g(t)G(t)^{a-1}[1-G(t)^{a}]^{b-1}\left\{1-[1-G(t)^{a}]^{b}\right\}^{-1}$$

and

$$H^{KwG}(t)=-b\log[1-G(t)^{a}].$$

Recently, Alizadeh *et al.*, (2015) proposed the Kumaraswamy Marshal-Olkin family of distributions by using the Marshal-Olkin (Marshall and Olkin 1997) cdf in that Kumaraswamy-*G* (Cordeiro and de Castro 2011) family and studied its many properties.

The main motivation behind the present article is to propose another family of continuous probability distribution that generalizes the Kumaraswamy-*G* (Cordeiro and de Castro 2011) family and the Marshall-Olkin Extended family (Marshall and Olkin 1997) by integrating the former as the base line distribution in the later. We call this new family the Marshall-Olkin Kumaraswamy-*G* $(MOKw-G)$ family of distribution which encompasses many known families of distributions and study some its general properties, parameter estimation and real life application in the present article.

The rest of this article is organized in seven sections. In section 2 the new family is defined along with its physical basis. Next section is devoted to presenting some important special cases of the family along with their shape and main reliability characteristics. In section 4 we discuss some important general results of the proposed family, while different methods of estimation of parameters are presented in section 5. In section 6 we present two examples of comparative data fitting. The article ends with a concluding discussion and remark in the final section followed by an appendix to derive asymptotic confidence bounds.

**2. New Marshall-Olkin Kumaraswamy-*G* $(MOKw-G)$ family of distributions**

We now propose a new extension of the $MO$ (Marshall and Olkin 1997) family by considering the cdf and pdf of $Kw-G$ (Cordeiro and de Castro 2011) distribution in (4) and (5) as the $f(t)$ and $F(t)$ respectively in the $MO$ formulation in (3) and call it $MOKw-G$ distribution. The pdf of $MOKw-G$ is given by

$$\begin{aligned} f^{MOKwG}(t) &= \alpha\, a b\, g(t) G(t)^{a-1}[1-G(t)^a]^{b-1} \big/ [\alpha + \overline{\alpha}\{1-[1-G(t)^a]^b\}]^2 \\ &= \alpha\, a b\, g(t) G(t)^{a-1}[1-G(t)^a]^{b-1} \big/ [(\alpha+\overline{\alpha}) - \overline{\alpha}[1-G(t)^a]^b]^2 \\ &= \alpha\, a b\, g(t) G(t)^{a-1}[1-G(t)^a]^{b-1} \big/ [1-\overline{\alpha}[1-G(t)^a]^b]^2 \qquad (6) \\ &\qquad , \; -\infty < t < \infty, \alpha > 0, a > 0, b > 0 \end{aligned}$$

Similarly using equation (4) in (2) we get the cdf and hence the sf of $MOKw-G$ respectively as

$$F^{MOKwG}(t) = 1-[1-G(t)^a]^b \big/ \alpha + \overline{\alpha}\left\{1-[1-G(t)^a]^b\right\} = 1-[1-G(t)^a]^b \big/ 1-\overline{\alpha}[1-G(t)^a]^b \qquad (7)$$

and $\overline{F}^{MOKwG}(t) = 1 - F^{MOKwG}(t) = [1-G(t)^a]^b(1-\overline{\alpha}) \big/ 1-\overline{\alpha}[1-G(t)^a]^b$

$$= \alpha[1-G(t)^a]^b \big/ 1-\overline{\alpha}[1-G(t)^a]^b \qquad (8)$$

hrf: $h^{MOKwG}(t) = \dfrac{f^{MOKwG}(t)}{\overline{F}^{MOKwG}(t)} = \dfrac{\alpha\, a b\, g(t) G(t)^{a-1}[1-G(t)^a]^{b-1} \big/ [1-\overline{\alpha}[1-G(t)^a]^b]^2}{\alpha\,[1-G(t)^a]^b \big/ 1-\overline{\alpha}[1-G(t)^a]^b}$

$$= \frac{a\,b\,g(t)G(t)^{a-1}[1-G(t)^{a}]^{-1}}{1-\overline{\alpha}[1-G(t)^{a}]^{b}} \tag{9}$$

rhrf : $r^{MOKwG}(t) = \dfrac{f^{MOKwG}(t)}{F^{MOKwG}(t)} = \dfrac{\alpha\,a\,b\,g(t)G(t)^{a-1}[1-G(t)^{a}]^{b-1} \big/ [1-\overline{\alpha}[1-G(t)^{a}]^{b}]^{2}}{1-[1-G(t)^{a}]^{b} \big/ 1-\overline{\alpha}[1-G(t)^{a}]^{b}}$

$$= \frac{\alpha\,a\,b\,g(t)G(t)^{a-1}[1-G(t)^{a}]^{b-1}[1-[1-G(t)^{a}]^{b}]^{-1}}{1-\overline{\alpha}[1-G(t)^{a}]^{b}}$$

chrf: $H^{MOKwG}(t) = -\log[\overline{F}^{MOKwG}(t)] = -\log\left\{\alpha[1-G(t)^{a}]^{b} \big/ 1-\overline{\alpha}[1-G(t)^{a}]^{b}\right\}$

For

(i) $\alpha = 1$, the pdf in (6) reduces to that of $Kw-G$ (Cordeiro and de Castro, 2011)

(ii) $a = b = 1$, the pdf in (6) reduces to that of $MO$ (Marshall and Olkin, 1997).

**2.1 Genesis of the distribution**

Let $T_i$ $(i = 1, 2, \cdots, N)$ be a sequence of *i. i. d.* random variables with survival function $[1-G(t)^{a}]^{b}$, and let $N$ be a positive integer random variable independent of the $T_i$ 's given by the pgf of a geometric distribution with parameter $\alpha$, say $\vartheta(x;\alpha) = \alpha\,x(1-\overline{\alpha}\,x)^{-1}$. Then the inverse of $\vartheta(x;\alpha)$ becomes $\vartheta^{-1}(x;\alpha) = \vartheta(x;\alpha^{-1})$.We can verify that equation (8) comes from $\overline{F}^{MOKwG}(t) = \vartheta[\{1-G(t)^{a}\}^{b};\alpha]$ for $0 < \alpha < 1$ and $\overline{F}^{MOKwG}(t) = \vartheta[\{1-G(t)^{a}\}^{b};\alpha^{-1}]$ for $\alpha > 1$. For both cases, equation (8) represents the survival function of the minimum and maximum of $T_i$ $(i = 1, 2, \cdots, N)$, where $N$ has probability generating function $\vartheta(x;.)$ with parameter $\alpha$ *and* $\alpha^{-1}$ respectively.

**2.2 Shape of the density and hazard functions**

Here we have plotted the pdf and hrf of the $MOKw-G$ for some choices of the distribution $G$ and parameter values to study the variety of shapes assumed by the family.

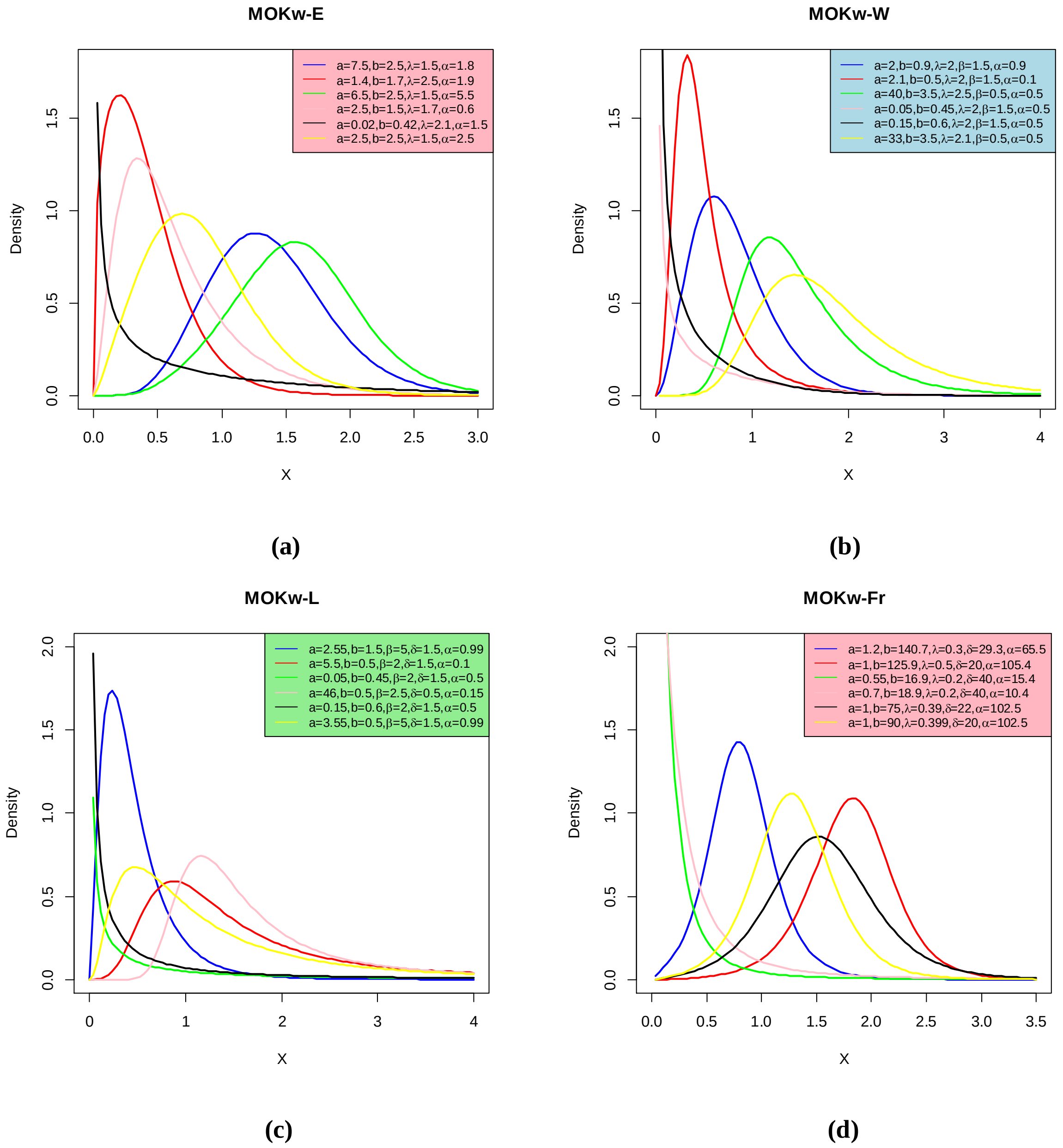

**Fig 1:** Density plots **(a)** $MOKw-E$ , **(b)** $MOKw-W$ , **(c)** $MOKw-L$ and **(d)** $MOKw-Fr$ distributions.

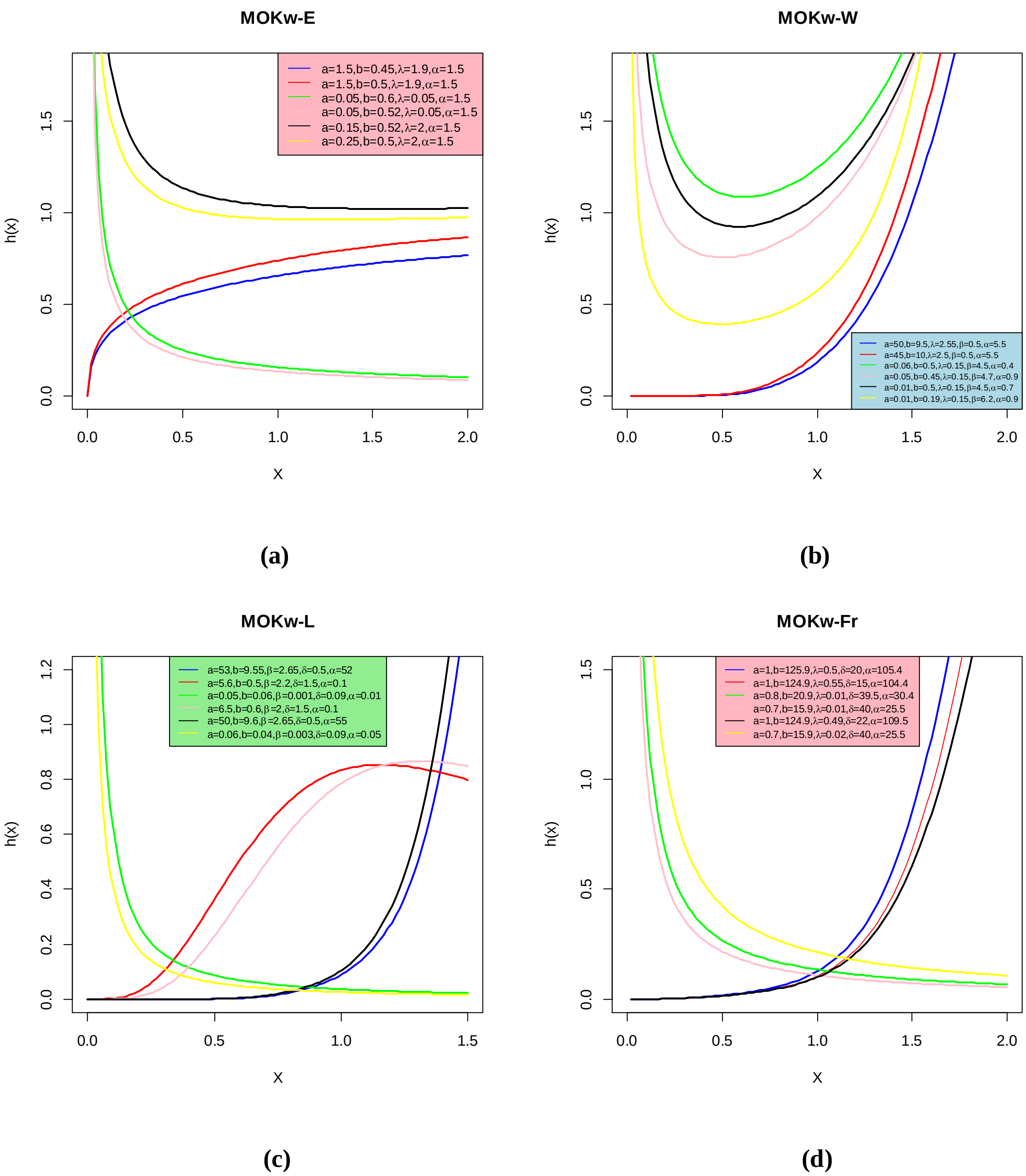

**Fig 2:** Hazard plots **(a)** $MOKw-E$ , **(b)** $MOKw-W$ , **(c)** $MOKw-L$ and **(d)** $MOKw-Fr$ distributions.

From the plots in figure 1 and 2 it can be seen that the family is very flexible and can offer many different types of shapes of density and hazard rate function including the bath tub shaped for hazard.

## 3. Special models

In this section we provide some special cases of the $MOKw-G$ family of distributions and list their main distributional characteristics.

### 3.1 The $MOKw-$exponential ($MOKw-E$) distribution

Let the base line distribution be exponential with parameter $\lambda>0,\ g(t)=\lambda e^{-\lambda t},\ t>0$ and $G(t)=1-e^{-\lambda t},\quad t>0$ then for the $MOKw-E$ model we get

pdf: $$f^{MOKwE}(t)=\alpha\, a b \lambda e^{-\lambda t}(1-e^{-\lambda t})^{a-1}[1-(1-e^{-\lambda t})^{a}]^{b-1}\big/[1-\bar{\alpha}[1-(1-e^{-\lambda t})^{a}]^{b}]^{2}$$

cdf : $$F^{MOKwE}(t)=1-[1-(1-e^{-\lambda t})^{a}]^{b}\big/1-\bar{\alpha}[1-(1-e^{-\lambda t})^{a}]^{b}$$

sf : $$\bar{F}^{MOKwE}(t)=\alpha[1-(1-e^{-\lambda t})^{a}]^{b}\big/1-\bar{\alpha}[1-(1-e^{-\lambda t})^{a}]^{b}$$

hrf : $$h^{MOKwE}(t)=\frac{a b \lambda e^{-\lambda t}(1-e^{-\lambda t})^{a-1}[1-(1-e^{-\lambda t})^{a}]^{-1}}{1-\bar{\alpha}[1-(1-e^{-\lambda t})^{a}]^{b}}$$

rhrf: : $$r^{MOKwE}(t)=\frac{\alpha\, a b \lambda e^{-\lambda t}(1-e^{-\lambda t})^{a-1}[1-(1-e^{-\lambda t})^{a}]^{b-1}\{1-[1-(1-e^{-\lambda t})^{a}]^{b}\}^{-1}}{1-\bar{\alpha}[1-(1-e^{-\lambda t})^{a}]^{b}}$$

chrf : $$H^{MOKwE}(t)=-\log[\bar{F}^{MOKwE}(t)]=-\log[\alpha[1-(1-e^{-\lambda t})^{a}]^{b}\big/1-\bar{\alpha}[1-(1-e^{-\lambda t})^{a}]^{b}]$$

### 3.2 The $MOKw-$Lomax ($MOKw-L$) distribution

Considering the Lomax distribution (Ghitany *et al.*, 2007) with pdf and cdf given by $g(t)=(\beta/\delta)[1+(t/\delta)]^{-(\beta+1)},\quad t>0$, and $G(t)=1-[1+(t/\delta)]^{-\beta},\quad \beta>0$ and $\delta>0$ the pdf of the $MOKw-L$ distribution is given by

pdf: $$f^{MOKwL}(t)=\frac{\alpha\, a b(\beta/\delta)[1+(t/\delta)]^{-(\beta+1)}[1-[1+(t/\delta)]^{-\beta}]^{a-1}[1-[1-[1+(t/\delta)]^{-\beta}]^{a}]^{b-1}}{[1-\bar{\alpha}[1-[1-[1+(t/\delta)]^{-\beta}]^{a}]^{b}]^{2}}$$

cdf: $$F^{MOKwL}(t)=\frac{1-[1-\{1-[1+(t/\delta)]^{-\beta}\}^{a}]^{b}}{1-\bar{\alpha}[1-\{1-[1+(t/\delta)]^{-\beta}\}^{a}]^{b}}$$

sf: $$\bar{F}^{MOKwL}(t)=\frac{\alpha[1-\{1-[1+(t/\delta)]^{-\beta}\}^{a}]^{b}}{1-\bar{\alpha}[1-\{1-[1+(t/\delta)]^{-\beta}\}^{a}]^{b}}$$

hrf: $$h^{MOKwL}(t)=\frac{a b(\beta/\delta)[1+(t/\delta)]^{-(\beta+1)}[1-[1+(t/\delta)]^{-\beta}]^{a-1}[1-\{1-[1+(t/\delta)]^{-\beta}\}^{a}]^{-1}}{1-\bar{\alpha}[1-\{1-[1+(t/\delta)]^{-\beta}\}^{a}]^{b}}$$

rhrf: $r^{MOKwL}(t)=$

$$\left\{\begin{matrix}\alpha\, a b(\beta/\delta)[1+(t/\delta)]^{-(\beta+1)}[1-[1+(t/\delta)]^{-\beta}]^{a-1} \\ [1-\{1-[1+(t/\delta)]^{-\beta}\}^{a}]^{b-1}[1-[1-\{1-[1+(t/\delta)]^{-\beta}\}^{a}]^{b}]^{-1}\end{matrix}\right\}\Bigg/1-\bar{\alpha}[1-\{1-[1+(t/\delta)]^{-\beta}\}^{a}]^{b}$$

chrf: $$H^{MOKwL}(t)=-\log[\bar{F}^{MOKwE}(t)]=-\log\left[\frac{\alpha[1-\{1-[1+(t/\delta)]^{-\beta}\}^{a}]^{b}}{1-\bar{\alpha}[1-\{1-[1+(t/\delta)]^{-\beta}\}^{a}]^{b}}\right]$$

### 3.3 The $MOKw-$Weibull ($MOKw-W$) distribution

Considering the Weibull distribution (Ghitany *et al.*, 2005, Zhang and Xie, 2007, and Caroni, 2010) with parameters $\lambda > 0$ and $\beta > 0$ having pdf and cdf $g(t) = \lambda\beta t^{\beta-1} e^{-\lambda t^{\beta}}$ and $G(t) = 1 - e^{-\lambda t^{\beta}}$ respectively we get the pdf of $MOKw - W$ distribution as

$$f^{MOKwW}(t) = \alpha ab\lambda\beta t^{\beta-1} e^{-\lambda t^{\beta}} [1 - e^{-\lambda t^{\beta}}]^{a-1} [1 - [1 - e^{-\lambda t^{\beta}}]^a]^{b-1} \Big/ [1 - \overline{\alpha}[1 - [1 - e^{-\lambda t^{\beta}}]^a]^b]^2$$

cdf: $F^{MOKwW}(t) = \{1 - [1 - [1 - e^{-\lambda t^{\beta}}]^a]^b\} / \{1 - \overline{\alpha}[1 - [1 - e^{-\lambda t^{\beta}}]^a]^b\}$

sf: $\overline{F}^{MOKwW}(t) = \{\alpha[1 - [1 - e^{-\lambda t^{\beta}}]^a]^b\} / \{1 - \overline{\alpha}[1 - [1 - e^{-\lambda t^{\beta}}]^a]^b\}$

hrf: $h^{MOKwW}(t) = \{ab\lambda\beta t^{\beta-1} e^{-\lambda t^{\beta}} [1 - e^{-\lambda t^{\beta}}]^{a-1} [1 - [1 - e^{-\lambda t^{\beta}}]^a]^{-1}\} / \{1 - \overline{\alpha}[1 - [1 - e^{-\lambda t^{\beta}}]^a]^b\}$

rhrf: $r^{MOKwW}(t) = \dfrac{\alpha ab\lambda\beta t^{\beta-1} e^{-\lambda t^{\beta}} [1 - e^{-\lambda t^{\beta}}]^{a-1} [1 - [1 - e^{-\lambda t^{\beta}}]^a]^{b-1} [1 - [1 - [1 - e^{-\lambda t^{\beta}}]^a]^b]^{-1}}{1 - \overline{\alpha}[1 - [1 - e^{-\lambda t^{\beta}}]^a]^b}$

chrf: $H^{MOKwW}(t) = -\log[\overline{F}^{MOKwW}(t)] = -\log\left[\dfrac{\alpha[1 - [1 - e^{-\lambda t^{\beta}}]^a]^b}{1 - \overline{\alpha}[1 - [1 - e^{-\lambda t^{\beta}}]^a]^b}\right]$

**3.4 The $MOKw -$ Frechet ( $MOKw - Fr$ ) distribution**

Suppose the base line distribution is the Frechet distribution (Krishna *et al.*, 2013) with pdf and cdf given by $g(t) = \lambda\delta^{\lambda} t^{-(\lambda+1)} e^{-(\delta/t)^{\lambda}}$ and $G(t) = e^{-(\delta/t)^{\lambda}}, t > 0$ respectively, and then the corresponding pdf of $MOKw - Fr$ distribution becomes pdf: $f^{MOKwFr}(t) = \dfrac{\alpha ab\lambda\delta^{\lambda} t^{-(\lambda+1)} e^{-(\delta/t)^{\lambda}} [e^{-(\delta/t)^{\lambda}}]^{a-1} [1 - [e^{-(\delta/t)^{\lambda}}]^a]^{b-1}}{[1 - \overline{\alpha}[1 - [e^{-(\delta/t)^{\lambda}}]^a]^b]^2}$ cdf: $F^{MOKwFr}(t) = \dfrac{1 - [1 - [e^{-(\delta/t)^{\lambda}}]^a]^b}{1 - \overline{\alpha}[1 - [e^{-(\delta/t)^{\lambda}}]^a]^b}$ sf: $\overline{F}^{MOKwFr}(t) = \dfrac{\alpha[1 - (e^{-(\delta/t)^{\lambda}})^a]^b}{1 - \overline{\alpha}[1 - [e^{-(\delta/t)^{\lambda}}]^a]^b}$ hrf: $h^{MOKwFr}(t) = \dfrac{ab\lambda\delta^{\lambda} t^{-(\lambda+1)} e^{-(\delta/t)^{\lambda}} (e^{-(\delta/t)^{\lambda}})^{a-1} [1 - (e^{-(\delta/t)^{\lambda}})^a]^{-1}}{1 - \overline{\alpha}[1 - [e^{-(\delta/t)^{\lambda}}]^a]^b}$ rhrf: $r^{MOKwFr}(t) = \dfrac{\alpha ab\lambda\,\delta^{\lambda} t^{-(\lambda+1)} e^{-(\delta/t)^{\lambda}} (e^{-(\delta/t)^{\lambda}})^{a-1} [1 - (e^{-(\delta/t)^{\lambda}})^a]^{b-1} [1 - [1 - [e^{-(\delta/t)^{\lambda}}]^a]^b]^{-1}}{1 - \overline{\alpha}[1 - [e^{-(\delta/t)^{\lambda}}]^a]^b}$ chrf: $H^{MOKwFr}(t) = -\log[\overline{F}^{MOKwFr}(t)] = -\log\left[\dfrac{\alpha[1 - (e^{-(\delta/t)^{\lambda}})^a]^b}{1 - \overline{\alpha}[1 - [e^{-(\delta/t)^{\lambda}}]^a]^b}\right]$

**3.5 The $MOKw -$ Gompertz ( $MOKw - Go$ ) distribution**

Next by taking the Gompertz distribution (Gieser *et al.*, 1998) with pdf and cdf $g(t) = \beta e^{\lambda t} e^{-\frac{\beta}{\lambda}(e^{\lambda t} - 1)}$ and $G(t) = 1 - e^{-\frac{\beta}{\lambda}(e^{\lambda t} - 1)}, \beta > 0, \lambda > 0, t > 0$ respectively, we get the pdf of $MOKw - Go$ distribution as pdf: $f^{MOKwGo}(t)=\dfrac{\alpha ab\beta e^{\lambda t} e^{-\frac{\beta}{\lambda}(e^{\lambda t}-1)}[1-e^{-\frac{\beta}{\lambda}(e^{\lambda t}-1)}]^{a-1}[1-\{1-e^{-\frac{\beta}{\lambda}(e^{\lambda t}-1)}\}^{a}]^{b-1}}{[1-\overline{\alpha}[1-\{1-e^{-\frac{\beta}{\lambda}(e^{\lambda t}-1)}\}^{a}]^{b}]^{2}}$ cdf: $F^{MOKwGo}(t)=\dfrac{1-[1-\{1-e^{-\frac{\beta}{\lambda}(e^{\lambda t}-1)}\}^{a}]^{b}}{1-\overline{\alpha}[1-\{1-e^{-\frac{\beta}{\lambda}(e^{\lambda t}-1)}\}^{a}]^{b}}$ sf: $\overline{F}^{MOKwGo}(t)=\dfrac{\alpha[1-\{1-e^{-\frac{\beta}{\lambda}(e^{\lambda t}-1)}\}^{a}]^{b}}{1-\overline{\alpha}[1-\{1-e^{-\frac{\beta}{\lambda}(e^{\lambda t}-1)}\}^{a}]^{b}}$ hrf: $h^{MOKwGo}(t)=\dfrac{ab\beta e^{\lambda t} e^{-\frac{\beta}{\lambda}(e^{\lambda t}-1)}[1-e^{-\frac{\beta}{\lambda}(e^{\lambda t}-1)}]^{a-1}[1-[1-e^{-\frac{\beta}{\lambda}(e^{\lambda t}-1)}]^{a}]^{-1}}{1-\overline{\alpha}[1-\{1-e^{-\frac{\beta}{\lambda}(e^{\lambda t}-1)}\}^{a}]^{b}}$

rhrf: $r^{MOKwGo}(t)$

$$=\frac{\alpha ab\beta e^{\lambda t} e^{-\frac{\beta}{\lambda}(e^{\lambda t}-1)}[1-e^{-\frac{\beta}{\lambda}(e^{\lambda t}-1)}]^{a-1}[1-\{1-e^{-\frac{\beta}{\lambda}(e^{\lambda t}-1)}\}^{a}]^{b-1}[1-\{1-[1-e^{-\frac{\beta}{\lambda}(e^{\lambda t}-1)}]^{a}\}^{b}]^{-1}}{1-\overline{\alpha}[1-\{1-e^{-\frac{\beta}{\lambda}(e^{\lambda t}-1)}\}^{a}]^{b}}$$

chrf: $H^{MOKwGo}(t)=-\log[\overline{F}^{MOKwGo}(t)]=-\log\left[\dfrac{\alpha[1-\{1-e^{-\frac{\beta}{\lambda}(e^{\lambda t}-1)}\}^{a}]^{b}}{1-\overline{\alpha}[1-\{1-e^{-\frac{\beta}{\lambda}(e^{\lambda t}-1)}\}^{a}]^{b}}\right]$

**3.6 The $MOKw-EW$ distribution**

The extended Weibull ($EW$) distributions of Gurvich *et al.* (1997) has the cdf $G(t)=1-\exp[-\delta E(t:\xi)]$, $t\in D\subseteq R_{+},\delta>0$ where $E(t:\xi)$ is a non-negative monotonically increasing function which depends on the parameter vector $\xi$**.** The corresponding pdf is $g(t)=\delta\exp[-\delta E(t:\xi)]e(t:\xi)$ where $e(t:\xi)$ is the derivative of $E(t:\xi)$.

Some important particular cases of $EW$ can be seen as follows:

(i) $E(t:\xi)=t$: exponential distribution.

(ii) $E(t:\xi)=t^{2}$: Rayleigh (Burr type-X) distribution.

(iii) $E(t:\xi)=\log(t/k)$: Pareto distribution

(iv) $E(t:\xi)=\beta^{-1}[\exp(\beta t)-1]$: Gompertz distribution.

Here we derive $MOKw-EW$ by considering $EW$ as the base line distribution as

pdf: $f^{MOKwEW}(t)$

$$=\frac{\alpha ab\delta\exp[-\delta E(t:\xi)]e(t:\xi)[1-\exp[-\delta E(t:\xi)]]^{a-1}[1-\{1-\exp[-\delta E(t:\xi)]\}^{a}]^{b-1}}{[1-\overline{\alpha}[1-\{1-\exp[-\delta E(t:\xi)]\}^{a}]^{b}]^{2}}$$

cdf: $F^{MOKwEW}(t)=\dfrac{1-[1-\{1-\exp[-\delta E(t:\xi)]\}^{a}]^{b}}{1-\overline{\alpha}[1-\{1-\exp[-\delta E(t:\xi)]\}^{a}]^{b}}$

sf: $\overline{F}^{MOKwEW}(t) = \dfrac{\alpha\,[1-\{1-\exp[-\delta\,E(t:\xi)]\}^{a}]^{b}}{1-\overline{\alpha}\,[1-\{1-\exp[-\delta\,E(t:\xi)]\}^{a}]^{b}}$

hrf: $h^{MOKwEW}(t) = \dfrac{a\,b\,\delta\exp[-\delta\,E(t:\xi)]e(t:\xi)[1-\exp\{-\delta\,E(t:\xi)\}]^{a-1}\;[1-\{1-\exp[-\delta\,E(t:\xi)]\}^{a}]^{-1}}{1-\overline{\alpha}\,[1-\{1-\exp[-\delta\,E(t:\xi)]\}^{a}]^{b}}$

rhrf: $r^{MOKwEW}(t)$

$$= \frac{\alpha\,a\,b\,\delta\exp[-\delta\,E(t:\xi)]e(t:\xi)[1-\exp\{-\delta\,E(t:\xi)\}]^{a-1}\;[1-\{1-\exp[-\delta\,E(t:\xi)]\}^{a}]^{b-1}}{1-\overline{\alpha}\,[1-\{1-\exp[-\delta\,E(t:\xi)]\}^{a}]^{b}}$$

$$\times[1-[1-\{1-\exp[-\delta\,E(t:\xi)]\}^{a}]^{b}]^{-1}$$

chrf: $H^{MOKwEW}(t) = -\log\left[\dfrac{\alpha\,[1-\{1-\exp[-\delta\,E(t:\xi)]\}^{a}]^{b}}{1-\overline{\alpha}\,[1-\{1-\exp[-\delta\,E(t:\xi)]\}^{a}]^{b}}\right]$

**3.7 The $MOKw-$Extended modified Weibull ($MOKw-EMW$) distribution**

The cdf and pdf of the modified Weibull (*MW*) distribution (Sarhan and Zaindin, 2013) is given respectively by

$G(t) = 1-\exp[-\sigma\,t-\beta\,t^{\gamma}]$ **,** $t>0,\ \gamma>0, \sigma, \beta \geq 0, \sigma+\beta>0$ and

$g(t) = (\sigma+\beta\,\gamma\,t^{\gamma-1})\exp[-\sigma\,t-\beta\,t^{\gamma}]$.

The corresponding pdf of $MOKw-EMW$ is given by

$f^{MOKwEMW}(t)$

$$= \frac{\alpha\,a\,b(\sigma+\beta\,\gamma\,t^{\gamma-1})\exp(-\sigma\,t-\beta\,t^{\gamma})[1-\exp(-\sigma\,t-\beta\,t^{\gamma})]^{a-1}[1-[1-\exp(-\sigma\,t-\beta\,t^{\gamma})]^{a}]^{b-1}}{[1-\overline{\alpha}\{1-[1-\exp(-\sigma\,t-\beta\,t^{\gamma})]^{a}\}^{b}]^{2}}$$ cdf:

$F^{MOKwEMW}(t) = \dfrac{1-\{1-[1-\exp(-\sigma\,t-\beta\,t^{\gamma})]^{a}\}^{b}}{1-\overline{\alpha}\{1-[1-\exp(-\sigma\,t-\beta\,t^{\gamma})]^{a}\}^{b}}$

sf: $\overline{F}^{MOKwEMW}(t) = \dfrac{\alpha\{1-[1-\exp(-\sigma\,t-\beta\,t^{\gamma})]^{a}\}^{b}}{1-\overline{\alpha}\{1-[1-\exp(-\sigma\,t-\beta\,t^{\gamma})]^{a}\}^{b}}$

hrf: $h^{MOKwEMW}(t)$

$$= \frac{a\,b(\sigma+\beta\,\gamma\,t^{\gamma-1})\exp(-\sigma\,t-\beta\,t^{\gamma})[1-\exp(-\sigma\,t-\beta\,t^{\gamma})]^{a-1}[1-[1-\exp(-\sigma\,t-\beta\,t^{\gamma})]^{a}]^{-1}}{1-\overline{\alpha}\{1-[1-\exp(-\sigma\,t-\beta\,t^{\gamma})]^{a}\}^{b}}$$

rhrf: $r^{MOKwEMW}(t)$

$$= \frac{\alpha\,a\,b(\sigma+\beta\,\gamma\,t^{\gamma-1})\exp(-\sigma\,t-\beta\,t^{\gamma})[1-\exp(-\sigma\,t-\beta\,t^{\gamma})]^{a-1}[1-[1-\exp(-\sigma\,t-\beta\,t^{\gamma})]^{a}]^{b-1}}{1-\overline{\alpha}\{1-[1-\exp(-\sigma\,t-\beta\,t^{\gamma})]^{a}\}^{b}}$$

$$\times[1-\{1-[1-\exp(-\sigma\,t-\beta\,t^{\gamma})]^{a}\}^{b}]^{-1}$$

chrf: $H^{MOKwEMW}(t) = -\log[\bar{F}^{MOKwEMW}(t)] = -\log\left[\frac{\alpha\{1-[1-\exp(-\sigma t-\beta t^{\gamma})]^{a}\}^{b}}{1-\bar{\alpha}\{1-[1-\exp(-\sigma t-\beta t^{\gamma})]^{a}\}^{b}}\right]$

**3.8 The $MOKw-$Power Log-normal ($MOKw-PLn$) distribution**

Power Log-normal distribution introduced by Nelson and Dognanksoy (1992) by extending the Log-normal distribution and its density and cdf are respectively given by,

$$g(t) = \frac{p}{t\sigma}\Phi\left(\frac{\mu-\ln(t)}{\sigma}\right)\left[\Phi\left(\frac{\mu-\ln(t)}{\sigma}\right)\right]^{p-1} \text{ and } G(t) = 1-\left[\Phi\left(\frac{\mu-\ln(t)}{\sigma}\right)\right]^{p},$$

$-\infty<\mu<\infty, \sigma>0, p>0, t>0$.

For $p=1$, it reduces to the usual Log-normal distribution.

With this as the base line distribution we get the pdf of $MOKw-PLn$ distribution as

$$f^{MOKwPLn}(t) =$$

$$\alpha a b\frac{p}{t\sigma}\Phi\left(\frac{\mu-\ln(t)}{\sigma}\right)\left[\Phi\left(\frac{\mu-\ln(t)}{\sigma}\right)\right]^{p-1}\left[1-\left\{\Phi\left(\frac{\mu-\ln(t)}{\sigma}\right)\right\}^{p}\right]^{a-1}$$

$$\left[1-\left[1-\left\{\Phi\left(\frac{\mu-\ln(t)}{\sigma}\right)\right\}^{p}\right]^{a}\right]^{b-1}\Bigg/\left[1-\bar{\alpha}\left[1-\left[1-\left\{\Phi\left(\frac{\mu-\ln(t)}{\sigma}\right)\right\}^{p}\right]^{a}\right]^{b}\right]^{2}$$

cdf: $$F^{MOKwPLn}(t) = \left[1-\left[1-\left[1-\left\{\Phi\left(\frac{\mu-\ln(t)}{\sigma}\right)\right\}^{p}\right]^{a}\right]^{b}\right]\Bigg/\left[1-\bar{\alpha}\left[1-\left[1-\left\{\Phi\left(\frac{\mu-\ln(t)}{\sigma}\right)\right\}^{p}\right]^{a}\right]^{b}\right]$$

sf: $\bar{F}^{MOKwPLn}(t)$

$$= \left[\alpha\left[1-\left[1-\left\{\Phi\left(\frac{\mu-\ln(t)}{\sigma}\right)\right\}^{p}\right]^{a}\right]^{b}\right]\Bigg/\left[1-\bar{\alpha}\left[1-\left[1-\left\{\Phi\left(\frac{\mu-\ln(t)}{\sigma}\right)\right\}^{p}\right]^{a}\right]^{b}\right]$$

hrf: $h^{MOKwPLn}(t)$

$$= \frac{a b\frac{p}{t\sigma}\Phi\left(\frac{\mu-\ln(t)}{\sigma}\right)\left[\Phi\left(\frac{\mu-\ln(t)}{\sigma}\right)\right]^{p-1}\left[1-\left\{\Phi\left(\frac{\mu-\ln(t)}{\sigma}\right)\right\}^{p}\right]^{a-1}\left[1-\left[1-\left\{\Phi\left(\frac{\mu-\ln(t)}{\sigma}\right)\right\}^{p}\right]^{a}\right]^{-1}}{1-\bar{\alpha}\left[1-\left[1-\left\{\Phi\left(\frac{\mu-\ln(t)}{\sigma}\right)\right\}^{p}\right]^{a}\right]^{b}}$$

rhrf: $r^{MOKwPLn}(t)$

$$= \frac{\frac{\alpha a b p}{t\sigma}\Phi\left(\frac{\mu-\ln(t)}{\sigma}\right)\left[\Phi\left(\frac{\mu-\ln(t)}{\sigma}\right)\right]^{p-1}\left[1-\left\{\Phi\left(\frac{\mu-\ln(t)}{\sigma}\right)\right\}^{p}\right]^{a-1}\left[1-\left[1-\left\{\Phi\left(\frac{\mu-\ln(t)}{\sigma}\right)\right\}^{p}\right]^{a}\right]^{b-1}}{1-\bar{\alpha}\left[1-\left[1-\left\{\Phi\left(\frac{\mu-\ln(t)}{\sigma}\right)\right\}^{p}\right]^{a}\right]^{b}}$$

$$\times\left[1-\left[1-\left[1-\left\{\Phi\left(\frac{\mu-\ln(t)}{\sigma}\right)\right\}^{p}\right]^{a}\right]^{b}\right]^{-1}$$

chrf: $H^{MOKwPLn}(t) = -\log\,[\bar{F}^{MOKwPLn}(t)]$

$$= -\log\left[\alpha\left[1-\left[1-\left\{\Phi\left(\frac{\mu-\ln(t)}{\sigma}\right)\right\}^{p}\right]^{a}\right]^{b}\Big/ 1-\bar{\alpha}\left[1-\left[1-\left\{\Phi\left(\frac{\mu-\ln(t)}{\sigma}\right)\right\}^{p}\right]^{a}\right]^{b}\right]$$

for $p=1$, it reduces to the $MOKw$-Log-normal distribution.

### 3.9 The $MOKw-$Extended Exponentiated Pareto ($MOKw-EEP$) distribution

The pdf and cdf of the exponentiated Pareto distribution, of Nadarajah (2005) are given respectively by

$g(t) = \gamma k \theta^{k} t^{-(k+1)}[1-(\theta/t)^{k}]^{\gamma-1}$ and $G(t) = [1-(\theta/t)^{k}]^{\gamma}$, $x > \theta$ ; $\theta, k, \gamma > 0$.

Thus the pdf of $MOKw-EEP$ distribution is given by

$$f^{MOKwEEP}(t) = \frac{\alpha a b \gamma k \theta^{k} t^{-(k+1)}\{1-(\theta/t)^{k}\}^{\gamma-1}[\{1-(\theta/t)^{k}\}^{\gamma}]^{a-1}[1-[\{1-(\theta/t)^{k}\}^{\gamma}]^{a}]^{b-1}}{[1-\bar{\alpha}[1-\{[1-(\theta/t)^{k}]^{\gamma}\}^{a}]^{b}]^{2}}$$

cdf: $F^{MOKwEEP}(t) = \dfrac{1-[1-\{[1-(\theta/t)^{k}]^{\gamma}\}^{a}]^{b}}{1-\bar{\alpha}[1-\{[1-(\theta/t)^{k}]^{\gamma}\}^{a}]^{b}}$ sf: $\bar{F}^{MOKwEEP}(t) = \dfrac{\alpha[1-\{[1-(\theta/t)^{k}]^{\gamma}\}^{a}]^{b}}{1-\bar{\alpha}[1-\{[1-(\theta/t)^{k}]^{\gamma}\}^{a}]^{b}}$

hrf: $h^{MOKwEEP}(t) = \dfrac{a b \gamma k \theta^{k} t^{-(k+1)}[1-(\theta/t)^{k}]^{\gamma-1}[1-(\theta/t)^{k}]^{\gamma(a-1)}[1-[1-(\theta/t)^{k}]^{\gamma a}]^{-1}}{1-\bar{\alpha}[1-[1-(\theta/t)^{k}]^{\gamma a}]^{b}}$

rhrf: $r^{MOKwEEP}(t) = \dfrac{\alpha a b \gamma k \theta^{k} t^{-(k+1)}[1-(\theta/t)^{k}]^{\gamma-1}[1-(\theta/t)^{k}]^{\gamma(a-1)}[1-[1-(\theta/t)^{k}]^{\gamma a}]^{b-1}}{1-\bar{\alpha}[1-[1-(\theta/t)^{k}]^{\gamma a}]^{b}}$

$$\times[1-[1-\{1-(\theta/t)^{k}\}^{\gamma a}]^{b}]^{-1}$$

chrf: $H^{MOKwEEP}(t) = -\log\left[\dfrac{\alpha[1-\{[1-(\theta/t)^{k}]^{\gamma}\}^{a}]^{b}}{1-\bar{\alpha}[1-\{[1-(\theta/t)^{k}]^{\gamma}\}^{a}]^{b}}\right]$

### 3.10 The $MOKw-$Extended Power ($MOKw-EP$) distribution

The cdf and pdf of the extended power distribution are $G(t)=(\theta t)^k$ and $g(t)=k\theta^k t^{k-1}$, $t\in(0,1/\theta)$ and $\theta>0$. The corresponding pdf of $MOKw-EP$ distribution is then given by

pdf: $$f^{MOKwEP}(t)=\frac{\alpha abk\theta^k t^{k-1}[(\theta t)^k]^{a-1}[1-[(\theta t)^k]^a]^{b-1}}{[1-\overline{\alpha}\{1-[(\theta t)^k]^a\}^b]^2}$$

cdf: $$F^{MOKwEP}(t)=\frac{1-[1-[(\theta t)^k]^a]^b}{1-\overline{\alpha}[1-[(\theta t)^k]^a]^b},$$ sf: $$\overline{F}^{MOKwEP}(t)=\frac{\alpha[1-[(\theta t)^k]^a]^b}{1-\overline{\alpha}[1-[(\theta t)^k]^a]^b}$$

rhf: $$h^{MOKwEP}(t)=\frac{abk\theta^k t^{k-1}[(\theta t)^k]^{a-1}[1-[(\theta t)^k]^a]^{-1}}{1-\overline{\alpha}[1-[(\theta t)^k]^a]^b}$$

rhrf: $$r^{MOKwEP}(t)=\frac{\alpha abk\theta^k t^{k-1}[(\theta t)^k]^{a-1}[1-[(\theta t)^k]^a]^{b-1}[1-\{1-[(\theta t)^k]^a\}^b]^{-1}}{1-\overline{\alpha}[1-[(\theta t)^k]^a]^b}$$

chrf: $$H^{MOKwEP}(t)=-\log\left[\frac{\alpha[1-[(\theta t)^k]^a]^b}{1-\overline{\alpha}[1-[(\theta t)^k]^a]^b}\right]$$

## 4. General results for the Marshall-Olkin's Kumaraswamy-*G* ($MOKw-G$) family of distributions

In this section we derive certain general results for the proposed $MOKw-G$ family following the methods described in Barreto-Souza *et al.* (2013), Cordeiro *et al.* (2014b), and Alizadeh *et al.* (2015).

### 4.1 Series expansions

Consider the series representation $(1-z)^{-k}=\sum_{j=0}^{\infty}\frac{\Gamma(k+j)}{\Gamma(k)\,j!}z^j$. (10)

This is valid for $|z|<1$ and $k>0$, where $\Gamma(.)$ is the gamma function. If $\alpha\in(0,1)$ using (10) in (6), we obtain

$$\begin{aligned}f^{MOKwG}(t)&=\alpha abg(t)G(t)^{a-1}[1-G(t)^a]^{b-1}[1-\overline{\alpha}[1-G(t)^a]^b]^{-2}\\&=\alpha abg(t)G(t)^{a-1}[1-G(t)^a]^{b-1}\sum_{j=0}^{\infty}(j+1)[\overline{\alpha}[1-G(t)^a]^b]^j\\&=\alpha abg(t)G(t)^{a-1}[1-G(t)^a]^{b-1}\sum_{j=0}^{\infty}(j+1)\,(1-\alpha)^j\{[1-G(t)^a]^b\}^j\\&=f^{KwG}(t)\sum_{j=0}^{\infty}\kappa_j[\overline{F}^{KwG}(t)]^j \qquad (11)\\&=\sum_{j=0}^{\infty}\kappa_j^{/}\frac{d}{dt}[\overline{F}^{KwG}(t)]^{j+1}\end{aligned}$$

Where $\kappa_j = \kappa_j(\alpha) = (j+1)\alpha(1-\alpha)^j ; \kappa'_j = \kappa'_j(\alpha) = -\alpha(1-\alpha)^j$ and

$\overline{F}^{KwG}(t) = 1 - F^{KwG}(t) = [1-G(t)^a]^b$ is the sf of $Kw-G$ distribution.

Again, $f^{MOKwG}(t) = \alpha\, a b\, g(t) G(t)^{a-1}[1-G(t)^a]^{b-1} \sum_{j=0}^{\infty}(j+1)\,(1-\alpha)^j [1-[1-\{[1-G(t)^a]^b\}^j]]$

$$= \alpha\, a b\, g(t) G(t)^{a-1}[1-G(t)^a]^{b-1} \sum_{j=0}^{\infty}(j+1)(1-\alpha)^j \sum_{k=0}^{j} \binom{j}{k} [-\{1-[1-G(t)^a]^b\}]^{j-k}$$

$$= \alpha\, a\, b\, g(t) G(t)^{a-1}[1-G(t)^a]^{b-1} \sum_{j=0}^{\infty}(j+1)(1-\alpha)^j (-1)^{j-k} \sum_{k=0}^{j} \binom{j}{k} [\{1-[1-G(t)^a]^b\}]^{j-k}$$

$$= f^{KwG}(t) \sum_{j=0}^{\infty} \sum_{k=0}^{j} \alpha(1-\alpha)^j (-1)^{j-k}(j+1)\binom{j}{k}[F^{KwG}(t)]^{j-k}$$

$$= f^{KwG}(t) \sum_{j=0}^{\infty} \sum_{k=0}^{j} \varphi_{j,k}[F^{KwG}(t)]^{j-k} \qquad (12)$$

$$= \sum_{j=0}^{\infty} \sum_{k=0}^{j} \frac{\varphi_{j,k}}{j-k+1} \frac{d}{dt}[F^{KwG}(t)]^{j-k+1}$$

$$= \sum_{j=0}^{\infty} \sum_{k=0}^{j} \varphi'_{j,k} \frac{d}{dt}[F^{KwG}(t)]^{j-k+1} \qquad (13)$$

Where $\varphi_{j,k} = \varphi_{j,k}(\alpha) = \alpha\,(1-\alpha)^j\,(-1)^{j-k}\,(j+1)\binom{j}{k} ; \varphi'_{j,k} = \varphi'_{j,k}(\alpha) = \alpha\,(1-\alpha)^j\,(-1)^{j-k}\binom{j+1}{k}$

and $F^{KwG}(t) = 1-[1-G(t)^a]^b$ is the cdf of $Kw-G$ (Cordeiro and de Castro 2011) distribution in equation (4). Similarly an expansion for the survival function of $MOKw-G$ [for $\alpha \in (0,1)$ ] can be derives as

$$\overline{F}^{MOKwG}(t) = \alpha[1-G(t)^a]^b \big/ 1-\overline{\alpha}[1-G(t)^a]^b$$

$$= \alpha[1-G(t)^a]^b \{1-\overline{\alpha}[1-G(t)^a]^b\}^{-1}$$

$$= \alpha[1-G(t)^a]^b \sum_{j=0}^{\infty} \frac{\Gamma(j+1)}{\Gamma(1)\, j!} \{\overline{\alpha}[1-G(t)^a]^b\}^j$$

$$= \alpha\, \overline{F}^{KwG}(t;a,b) \sum_{j=0}^{\infty} (1-\alpha)^j \{\overline{F}^{KwG}(t;a,b)\}^j$$

$$= -\sum_{j=0}^{\infty} \kappa'_j [\overline{F}^{KwG}(t;a,b)]^{j+1}$$

$$= -\sum_{j=0}^{\infty} \kappa'_j\, \overline{F}^{KwG}(t;a,b(j+1))$$

The density function (6) can also be expressed as

$$f^{MOKwG}(t)=\frac{a\,b\,g(t)G(t)^{a-1}[1-G(t)^{a}]^{b-1}}{\alpha\left[1-\frac{(\alpha-1)[1-[1-G(t)^{a}]^{b}]}{\alpha}\right]^{2}}$$

Hence for $\alpha>1$ using (10) we get

$$=\frac{a\,b\,g(t)G(t)^{a-1}[1-G(t)^{a}]^{b-1}}{\alpha}\left[1-\frac{(\alpha-1)[1-[1-G(t)^{a}]^{b}]}{\alpha}\right]^{-2}$$

$$=a\,b\,g(t)G(t)^{a-1}[1-G(t)^{a}]^{b-1}\sum_{j=0}^{\infty}(j+1)\left(\frac{\alpha-1}{\alpha}\right)^{j}\left\{1-[1-G(t)^{a}]^{b}\right\}^{j}$$

$$=\sum_{j=0}^{\infty}\left(\frac{\alpha-1}{\alpha}\right)^{j}(j+1)a\,b\,g(t)G(t)^{a-1}[1-G(t)^{a}]^{b-1}\left\{1-[1-G(t)^{a}]^{b}\right\}^{j}$$

$$=f^{KwG}(t)\sum_{j=0}^{\infty}\eta_{j}\left\{F^{KwG}(t)\right\}^{j} \tag{14}$$

$$=\sum_{j=0}^{\infty}\eta'_{j}\frac{d}{dt}\left\{F^{KwG}(t)\right\}^{j+1} \tag{15}$$

where $\eta_j=\eta_j(\alpha)=\frac{j+1}{\alpha}\left(1-\frac{1}{\alpha}\right)^{j}$ and $\eta'_j=\eta'_j(\alpha)=\frac{1}{\alpha}\left(1-\frac{1}{\alpha}\right)^{j}$

For $\alpha>1$, the survival function of $MOKw-G$ can be expressed as

$$\overline{F}^{MOKwG}(t)=\frac{[1-G(t)^{a}]^{b}}{1-(\alpha-1)[1-[1-G(t)^{a}]^{b}]/\alpha}$$

On using (10) we get

$$=[1-G(t)^{a}]^{b}\left[1-\frac{(\alpha-1)[1-[1-G(t)^{a}]^{b}]}{\alpha}\right]^{-1}$$

$$=[1-G(t)^{a}]^{b}\sum_{j=0}^{\infty}\frac{\Gamma(j+1)}{\Gamma(1)\;j!}\frac{(\alpha-1)^{j}}{\alpha^{j}}[1-[1-G(t)^{a}]^{b}]^{j}$$

$$=\overline{F}^{KwG}(t;a,b)\sum_{j=0}^{\infty}(-1)^{j}\left(\frac{\overline{\alpha}}{\alpha}\right)^{j}[F^{KwG}(t;a,b)]^{j}$$

$$=\overline{F}^{KwG}(t;a,b)\sum_{j=0}^{\infty}C'_{j}[F^{KwG}(t;a,b)]^{j}\ \text{ where } C'_{j}=(\alpha/j+1)\,\eta_{j}\ .$$

**4.2 Order statistics**

Suppose $T_1, T_2, ...T_n$ is a random sample from any $MOKw-G$ distribution. Let $T_{i:n}$ denote the $i^{th}$ order statistics. The pdf of $T_{i:n}$ can be expressed as

$$\begin{aligned} f_{i:n}(t) &= \frac{n!}{(i-1)!\,(n-i)!} f^{MOKwG}(t)[1-\overline{F}^{MOKwG}(t)]^{i-1}\,\overline{F}^{MOKwG}(t)^{n-i} \\ &= \frac{n!}{(i-1)!\,(n-i)!} f^{MOKwG}(t)\,\overline{F}^{MOKwG}(t)^{n-i} \sum_{p=0}^{i-1} \binom{i-1}{p} [-\overline{F}^{MOKwG}(t)]^{p} \\ &= \frac{n!}{(i-1)!\,(n-i)!} f^{MOKwG}(t) \sum_{p=0}^{i-1} (-1)^{p} \binom{i-1}{p} \overline{F}^{MOKwG}(t)^{n+p-i} \end{aligned}$$

Now using the general expansion of the $MOKw-G$ distribution pdf and sf we get the pdf of the $i^{th}$ order statistics for of the $MOKw-G$ for $\alpha \in (0,1)$ as

$f_{i:n}(t) =$

$$\frac{n!}{(i-1)!\,(n-i)!} \left\{ f^{KwG}(t) \sum_{j=0}^{\infty} \kappa_j [\overline{F}^{KwG}(t)]^{j} \right\} \sum_{p=0}^{i-1} (-1)^{p+1} \binom{i-1}{p} \left\{ \sum_{q=0}^{\infty} \kappa'_q [\overline{F}^{KwG}(t;a,b)]^{q+n+p+1-i} \right\}$$

where $\kappa_j$ and $\kappa'_q$ are defined earlier.

$$= \frac{n!}{(i-1)!\,(n-i)!} \sum_{p=o}^{i-1} \binom{i-1}{p} (-1)^{p+1} f^{KwG}(t;a,b) \sum_{j=0}^{\infty} \sum_{q=0}^{\infty} \kappa_j \kappa'_q [\overline{F}^{KwG}(t;a,b)]^{j+q+n+p+1-i}$$

$$= f^{KwG}(t;a,b) \sum_{j,q=0}^{\infty} X_{j,q} [\overline{F}^{KwG}(t;a,b)]^{j+q+n+p+1-i} \qquad (16)$$

$$= -\sum_{j,q=0}^{\infty} [X_{j,q}/(j+q+n+p-i+2)] \frac{d}{dt} [\overline{F}^{KwG}(t;a,b)]^{j+q+n+p-i+2}$$

$$= \sum_{j,q=0}^{\infty} X'_{j,q} \frac{d}{dt} [\overline{F}^{KwG}(t;a,b)]^{j+q+n+p-i+2)}$$

$$= \sum_{j,q=0}^{\infty} X'_{j,q} \frac{d}{dt} [\overline{F}^{KwG}(t;a,b(j+q+n+p-i+2))]$$

$$= \sum_{j,q=0}^{\infty} X'_{j,q} f^{KwG}(t;a,b(j+q+n+p-i+2)) \qquad (17)$$

Where $X_{j,q} = n\kappa_j \kappa'_q \binom{n-1}{i-1} \sum_{p=0}^{i-1} \binom{i-1}{p} (-1)^{p+1}$ and $X'_{j,q} = -X_{j,q}/(j+q+n+p-i+2)$.

Again using the general expansion of the pdf and sf of $MOKw-G$ distribution we get the pdf of the $i^{th}$ order statistics for of the $MOKw-G$ for $\alpha > 1$ as

$$f_{i:n}(t)=\frac{n!}{(i-1)!\,(n-i)!}\left\{f^{KwG}(t;a,b)\sum_{j=0}^{\infty}\eta_j\left\{F^{KwG}(t;a,b)\right\}^j\right\}\sum_{p=0}^{i-1}(-1)^p\binom{i-1}{p}$$

$$\left\{[\bar{F}^{KwG}(t;a,b)]\sum_{k=0}^{\infty}C'_k\,[F^{KwG}(t;a,b)]^k\right\}^{n+p-i}$$

$$=\frac{n!}{(i-1)!\,(n-i)!}\left\{f^{KwG}(t;a,b)\sum_{j=0}^{\infty}\eta_j\left\{F^{KwG}(t;a,b)\right\}^j\right\}\sum_{p=0}^{i-1}(-1)^p\binom{i-1}{p}$$

$$[\bar{F}^{KwG}(t;a,b)]^{(n+p-i)}\left[\sum_{k=0}^{\infty}C'_k\,[F^{KwG}(t;a,b)]^k\right]^{n+p-i}$$

Where $\eta_j$ and $C'_k$ defined in section 4.1

Now, $\left[\sum_{k=0}^{\infty}C'_k\,[F^{KwG}(t;a,b)]^k\right]^{n+p-i}=\sum_{k=0}^{\infty}d_{n+p-i,k}\,[F^{KwG}(t;a,b)]^k$ (Nadarajah et. al 2015)

Where $d_{n+p-i,k}=\frac{1}{k\,C'_0}\sum_{h=1}^{k}[h(n+p-i-1)-k]C'_h\,d_{n+p-i,k-h}$

The density function of the $i^{th}$ order statistics of $MOKw-G$ distribution can be expressed as

$$f_{i:n}(t)=\frac{n!}{(i-1)!\,(n-i)!}\left\{f^{KwG}(t;a,b)\sum_{j=0}^{\infty}\eta_j\left\{F^{KwG}(t;a,b)\right\}^j\right\}\sum_{p=0}^{i-1}(-1)^p\binom{i-1}{p}$$

$$[\bar{F}^{KwG}(t;a,b)]^{(n+p-i)}\sum_{k=0}^{\infty}d_{n+p-i,k}\,[F^{KwG}(t;a,b)]^k$$

$$=\frac{n!}{(i-1)!\,(n-i)!}\sum_{p=0}^{i-1}(-1)^p\binom{i-1}{p}f^{KwG}(t;a,b)[\bar{F}^{KwG}(t;a,b)]^{(n+p-i)}$$

$$\sum_{j,k=0}^{\infty}\eta_j\,d_{n+p-i,k}\left\{F^{KwG}(t;a,b)\right\}^{j+k}$$

$$=f^{KwG}(t;a,b)\ [\bar{F}^{KwG}(t;a,b)]^{(n+p-i)}\sum_{j,k=0}^{\infty}\lambda_{j,k}\left\{F^{KwG}(t;a,b)\right\}^{j+k} \qquad (18)$$

Where $\lambda_{j,k}=n\,\eta_j\,d_{n+p-i,k}\binom{n-1}{i-1}\sum_{p=0}^{i-1}(-1)^p\binom{i-1}{p}$

**Remark I**: Equations (11)-(18) reveal that the density functions of the $MOKw-G$ distribution and that of its order statistics can be expressed as a product of the baseline density $f(x)$ with an infinite power series of $G(.)$ and also as a mixture of exponentiated-$G$ distributions under Lehman alternatives. These results play important role and may be used to obtain explicit expressions for the moments and moment generating function (mgf) of the $MOKw-G$ distribution and of its order statistics in a general framework and also for special models using the corresponding results of exponentiated - $G$ distributions.

### 4.3 Probability weighted moments

The probability weighted moments (PWMs), first proposed by Greenwood *et al.* (1979), are expectations of certain functions of a random variable whose mean exists. The $(p,q,r)^{th}$ PWM of $T$ is having a base line cdf $F(t)$ defined by

$$\Gamma_{p,q,r} = \int_{-\infty}^{\infty} t^p [F(t)]^q [1-F(t)]^r f(t)\,dt$$

From equations (11), (12) and (14) the $s^{th}$ moment of $T$ for $\alpha \in (0,1)$ and $\alpha > 1$, can be written either as

$$E(T^s) = \int_{-\infty}^{+\infty} t^s\, a b\, g(t) G(t)^{a-1} [1-G(t)^a]^{b-1} \sum_{j=0}^{\infty} \kappa_j [1-G(t)^a]^{bj}\, dt$$

$$= \sum_{j=0}^{\infty} \kappa_j \int_{-\infty}^{+\infty} t^s\, [1-G(t)^a]^{bj}\, a b\, g(t) G(t)^{a-1} [1-G(t)^a]^{b-1}\, dt$$

$$= \sum_{j=0}^{\infty} \kappa_j\, \Gamma_{s,\,0,\,j} \quad \text{or}$$

$$E(T^s) = \int_{-\infty}^{+\infty} t^s\, a b\, g(t) G(t)^{a-1} [1-G(t)^a]^{b-1} \sum_{j=0}^{\infty} \sum_{k=0}^{j} \varphi_{j,k} [1-[1-G(t)^a]^b]^{j-k}\, dt$$

$$= \sum_{j=0}^{\infty} \sum_{k=0}^{j} \varphi_{j,k} \int_{-\infty}^{+\infty} t^s [1-[1-G(t)^a]^b]^{j-k}\, a b\, g(t) G(t)^{a-1} [1-G(t)^a]^{b-1}\, dt$$

$$= \sum_{j=0}^{\infty} \sum_{k=0}^{j} \varphi_{j,k}\, \Gamma_{s,\,j-k,\,0} \quad \text{or as}$$

$$E(T^s) = \int_{-\infty}^{+\infty} t^s\, a b\, g(t)\, G(t)^{a-1} [1-G(t)^a]^{b-1} \sum_{j=0}^{\infty} \eta_j [1-[1-G(t)^a]^b]^j\, dt$$

$$= \sum_{j=0}^{\infty} \eta_j \int_{-\infty}^{+\infty} t^s [1-[1-G(t)^a]^b]^j\, a b\, g(t) G(t)^{a-1} [1-G(t)^a]^{b-1}\, dt$$

$$= \sum_{j=0}^{\infty} \eta_j\, \Gamma_{s,\,j,\,0}$$

Where $\Gamma_{p,q,r} = \int_{-\infty}^{\infty} t^p\, \{1-[1-G(t)^a]^b\}^q \{[1-G(t)^a]^b\}^r [a b\, g(t)\, G(t)^{a-1} [1-G(t)^a]^{b-1}]\, dt$

is the PWM of $MOKw-G$ the for the baseline distribution $G$.

The PWM can generally be used for estimating parameters quantiles of generalized distributions. These moments have low variance and don't possess severe biases, and they compare favourably with estimators obtained by maximum likelihood (Alizadeh *et al*., 2015).

Proceeding as above we can derive $s^{th}$ moment of the $i^{th}$ order statistic $T_{i:n}$, in a random sample of size *n* from $MOKw-G$ for $\alpha \in (0,1)$ and $\alpha > 1$, on using equations (16) and (18), as

$E(T^{s}{}_{i,n}) = \sum_{j,q=0}^{\infty} X_{j,q}\,\Gamma_{s,0,\ j+q+n+p-i}$ and $E(T^{s}{}_{i,n}) = \sum_{j,k=0}^{\infty} \lambda_{j,k}\,\Gamma_{s,\ j+k,\ n+p-i}$ respectively. Where the quantities $\varphi_{j,k}, \kappa_j, \eta_j, X_{j,q}$ and $\lambda_{j,k}$ are defined in section 4.1 and 4.2.

**4.4 Moment generating function**

The moment generating function of $MOKw-G$ family can be easily expressed in terms of those of the exponentiated $Kw-G$ (Cordeiro and de Castro, 2011) distribution using the results of section 4.1. For example using equation (15) it can be seen that

$$M_T(s) = E[e^{st}] = \int_{-\infty}^{\infty} e^{st} f(t)\,dt = \int_{-\infty}^{\infty} e^{st} \sum_{j=0}^{\infty} \eta_j^{/} \frac{d}{dt}\left\{F^{KwG}(t)\right\}^{j+1} dt$$

$$= \sum_{j=0}^{\infty} \eta_j^{/} \int_{-\infty}^{\infty} e^{st} \frac{d}{dt}\left\{F^{KwG}(t)\right\}^{j+1} dt = \sum_{j=0}^{\infty} \eta_j^{/} M_X(s)$$

Where *X* follows exponentiated $Kw-G$ (Cordeiro and de Castro, 2011) distribution.

**4.5 Renyi Entropy**

The entropy of a random variable is a measure of uncertainty. The Rényi entropy is defined as

$$I_R(\delta) = (1-\delta)^{-1} \log\left(\int_{-\infty}^{\infty} f(t)^{\delta}\, dt\right)$$

where $\delta > 0$ and $\delta \neq 1$ For furthers details, see Song (2001). For $\alpha \in (0,1)$ using expansion (10) in (6) we can write

$$f(t)^{\delta} = \left[\alpha\, a\, b\, g(t) G(t)^{a-1}[1-G(t)^{a}]^{b-1} \big/ [1-\overline{\alpha}\{1-G(t)^{a}\}^{b}]^{2}\right]^{\delta}$$

$$= \frac{\alpha^{\delta}[a\, b\, g(t) G(t)^{a-1}[1-G(t)^{a}]^{b-1}]^{\delta}}{\Gamma(2\delta)} \sum_{j=0}^{\infty} (1-\alpha)^{j}\, \Gamma(2\delta + j) \frac{[\{1-G(t)^{a}\}^{b}]^{j}}{j!}$$

Thus for $\alpha \in (0,1)$ the Rényi entropy of $T$ can be obtained as

$$I_R(\delta) = (1-\delta)^{-1} \log\left(\sum_{j=0}^{\infty} q_j \int_{-\infty}^{\infty} \left[a\, b\, g(t) G(t)^{a-1}[1-G(t)^{a}]^{b-1}\right]^{\delta} [\{1-G(t)^{a}\}^{b}]^{j}\, dt\right)$$

Where $q_j = q_j(\alpha) = \dfrac{\alpha^{\delta}(1-\alpha)^{j}\,\Gamma(2\delta + j)}{\Gamma(2\delta)\, j!}$

The density function (6) can be expressed as

$$f(t)^{MOKwG} = a\,b\,g(t)G(t)^{a-1}[1-G(t)^{a}]^{b-1} \Big/ \alpha\left[1-\frac{(\alpha-1)[1-[1-G(t)^{a}]^{b}]}{\alpha}\right]^{2}$$

Therefore

$$f(t)^{\delta} = \left\{a\,b\,g(t)G(t)^{a-1}[1-G(t)^{a}]^{b-1} \Big/ \alpha\left[1-\frac{(\alpha-1)[1-\{1-G(t)^{a}\}^{b}]}{\alpha}\right]^{2}\right\}^{\delta}$$

For $\alpha > 1,$ using expansion (10) we get,

$$f(t)^{\delta} = \frac{[a\,b\,g(t)G(t)^{a-1}[1-G(t)^{a}]^{b-1}]^{\delta}}{\alpha^{\delta+j}\,\Gamma(2\delta)}\sum_{j=0}^{\infty}(\alpha-1)^{j}\,\Gamma(2\delta+j)\,\frac{\{1-[1-G(t)^{a}]^{b}\}^{j}}{j!}$$

Thus for $\alpha > 1,$ the Rényi entropy of $T$ also can be obtained as

$$I_R(\delta) = (1-\delta)^{-1}\log\left(\sum_{j=0}^{\infty} r_j \int_{-\infty}^{\infty}[\,a\,b\,g(t)G(t)^{a-1}[1-G(t)^{a}]^{b-1}\,]^{\delta}[1-[1-G(t)^{a}]^{b}\,]^{j}\,dt\right)$$

Where $r_j = r_j(\alpha) = \dfrac{(\alpha-1)^{j}\,\Gamma(2\delta+j)}{\alpha^{\delta+j}\,\Gamma(2\delta)\,j!}$

**4.6 Quantile function and random sample generation**

We shall now present a formula for generating $MOKw-G$ random variable by using inversion method by inverting the cdf or the survival function.

$$F^{MO}(t) = \frac{F(t)}{\alpha+\bar{\alpha}F(t)} \Rightarrow F(t)\,[1-\bar{\alpha}\,F^{MO}(t)] = \alpha\,F^{MO}(t)$$

$$\Rightarrow F(t) = \frac{\alpha\,F^{MO}(t)}{1-\bar{\alpha}\,F^{MO}(t)} \;\Rightarrow t = F^{-1}\left[\frac{\alpha\,F^{MO}(t)}{1-\bar{\alpha}\,F^{MO}(t)}\right] \tag{19}$$

Now for the $Kw-G$ we have

$\bar{F}(t) = [1-G(t)^{a}]^{b}$. Therefore

$$\log\bar{F}(t) = b\log[1-G(t)^{a}] \Rightarrow 1-G(t)^{a} = \bar{F}(t)^{1/b} \Rightarrow a\log G(t) = \log\left[1-\bar{F}(t)^{1/b}\right]$$

$$\Rightarrow G(t) = [1-\{1-F(t)\}^{1/b}]^{1/a} \Rightarrow t = G^{-1}[1-\{1-F(t)\}^{1/b}]^{1/a} \tag{20}$$

Now combining (19) and (20) we get for $MOKw-G$

$$t = G^{-1}\left[1-\left\{1-\frac{\alpha\,u}{1-\bar{\alpha}\,u}\right\}^{1/b}\right]^{1/a} \tag{21}$$

Where $u \sim U(0,1)$

The $p^{\text{th}}$ Quantile $t_p$ for $MOKw-G$ can be easily obtained from (21) as

$$t_p = G^{-1}\left[1-\left\{1-\frac{\alpha\, p}{1-\bar{\alpha}\, p}\right\}^{1/b}\right]^{1/a}$$

For example, let the $G$ be exponential distribution with parameter $\lambda > 0$, having pdf and cdf as $g(t:\lambda) = \lambda e^{-\lambda t}, t > 0$ and $G(t:\lambda) = 1 - e^{-\lambda t}$, respectively. Then the $p^{th}$ quantile is obtained as $-(1/\lambda)\log[1-p]$. Therefore, the $p^{th}$ quantile $t_p$, of $MOKw-E$ is given by

$$t_p = -\frac{1}{\lambda}\log\left[1-\left\{1-\left(1-\frac{\alpha\, p}{1-\bar{\alpha}\, p}\right)^{1/b}\right\}^{1/a}\right]$$

### 4.7 Asymptotes and shapes

Here we investigate the asymptotic shapes of the proposed family following the method followed in Alizadeh *et al.*, (2015).

**Proposition 1**.The asymptotes of equations (6), (7) and (9) as $t \to 0$ are given by

$$f(t) \sim \frac{a b\, g(t) G(t)^{a-1}}{\alpha} \quad \text{as } G(t) \to 0$$

$$F(t) \sim 0 \quad \text{as } G(t) \to 0$$

$$h(t) \sim \frac{a b\, g(t) G(t)^{a-1}}{\alpha} \quad \text{as } G(t) \to 0$$

**Proposition 2**.The asymptotes of equations (6), (7) and (9) as $t \to \infty$ are given by

$$f(t) \sim \alpha a b\, g(t)[1-G(t)^a]^{b-1} \quad \text{as } t \to \infty$$

$$F(t) \sim 1-[1-G(t)^a]^b \quad \text{as } t \to \infty$$

$$h(t) \sim a\, b\, g(t)[1-G(t)^a]^{-1} \quad \text{as } t \to \infty$$

The shapes of the density and hazard rate functions can be described analytically. The critical points of the $MOKw-G$ density function are the roots of the equation:

$$\frac{d \log[f(t)]}{dt} = \frac{g'(t)}{g(t)} + (a-1)\frac{g(t)}{G(t)} + a(1-b)\frac{g(t)G(t)^{a-1}}{1-G(t)^a} - 2\,\frac{\bar{\alpha}\, a b\, g(t) G(t)^{a-1}[1-G(t)^a]^{b-1}}{1-\bar{\alpha}[1-G(t)^a]^b} \qquad (22)$$

The number of roots of (22) may be more than one. In particular, if $t = t_0$ is a root of (22) then it corresponds to a local maximum, a local minimum or a point of inflexion depending on whether $\psi(t_0) < 0, \psi(t_0) > 0 \ or \ \psi(t_0) = 0$ where

$$\psi(t) = \frac{d^2}{dt^2}\log[f(t)]$$

$$\psi(t) = \frac{g(t)g''(t) - [g'(t)]^2}{g(t)^2} + (a-1)\,\frac{G(t)g'(t) - g(t)^2}{G(t)}$$

$$+a(1-b)[\,\frac{g'(t)G(t)^{a-1}}{1-G(t)^a} + \frac{(a-1)\,g(t)^2\,G(t)^{a-2}}{1-G(t)^a} + \frac{a\,g(t)^2\,G(t)^{2a-2}}{[1-G(t)^a]^2}]$$

$$-\frac{2\overline{\alpha}\,ab[g'(t)G(t)^{a-1}[1-G(t)^a]^{b-1}]}{1-\overline{\alpha}[1-G(t)^a]^b} - \frac{2\overline{\alpha}\,ab[(a-1)g(t)^2\,G(t)^{a-2}[1-G(t)^a]^{b-1}]}{1-\overline{\alpha}[1-G(t)^a]^b}$$

$$+\frac{2\overline{\alpha}\,ab[a(b-1)\,g(t)^2\,G(t)^{2a-2}[1-G(t)^a]^{b-2}]}{1-\overline{\alpha}[1-G(t)^a]^b} + \left[\frac{\overline{\alpha}\,ab\,g(t)G(t)^{a-1}[1-G(t)^a]^{b-1}}{1-\overline{\alpha}[1-G(t)^a]^b}\right]^2$$

$$= \frac{g(t)g''(t) - [g'(t)]^2}{g(t)^2} + (a-1)\frac{G(t)g'(t) - g(t)^2}{G(t)}$$

$$+a(1-b)[\,\frac{g'(t)G(t)^{a-1}}{1-G(t)^a} + \frac{(a-1)\,g(t)^2\,G(t)^{a-2}}{1-G(t)^a} + \frac{a\,g(t)^2\,G(t)^{2a-2}}{[1-G(t)^a]^2}] - \frac{2\overline{\alpha}\,g'(t)}{g(t)}[1-G(t)^a]^b h^{MOKwG}(t)$$

$$-2\overline{\alpha}(a-1)\,g(t)G(t)^{-1}[1-G(t)^a]^b\,h^{MOKwG}(t) + 2\overline{\alpha}\,a(b-1)\,g(t)G(t)^{a-1}(t)\,\{1-G(t)^a\}^{b-1}h^{MOKwG}(t)$$

$$+\{[1-G(t)^a]^b\overline{\alpha}\,h^{MOKwG}(t)\}^2$$

The critical points of $h(t)$ are the roots of the equation

$$\frac{d\log[h(t)]}{dt} = \frac{g'(t)}{g(t)} + (a-1)\frac{g(t)}{G(t)} + \frac{a\,g(t)G(t)^{a-1}}{1-G(t)^a} - \frac{\overline{\alpha}\,ab\,g(t)G(t)^{a-1}[1-G(t)^a]^{b-1}}{1-\overline{\alpha}[1-G(t)^a]^b} \qquad (23)$$

The number of roots of (23) may be more than one. In particular, ff $t = t_0$ is a root of (23) then it corresponds to a local maximum, a local minimum or a point of inflexion depending on whether $\omega(t_0) < 0, \ \omega(t_0) > 0 \ or \ \omega(t_0) = 0$ where $\omega(t) = \frac{d^2}{dt^2}\log[h(t)]$

$$\omega(t) = \frac{g(t)g''(t) - [g'(t)]^2}{g(t)^2} + (a-1)\,\frac{G(t)g'(t) - g(t)^2}{G(t)}$$

$$+a\,[\,\frac{g'(t)G(t)^{a-1}}{1-G(t)^a} + \frac{(a-1)\,g(t)^2\,G(t)^{a-2}}{1-G(t)^a} + \frac{a\,g(t)^2\,G(t)^{2a-2}}{[1-G(t)^a]^2}]$$

$$-\frac{\overline{\alpha}\,ab[g'(t)G(t)^{a-1}[1-G(t)^a]^{b-1}]}{1-\overline{\alpha}[1-G(t)^a]^b} - \frac{\overline{\alpha}\,ab[(a-1)g(t)^2\,G(t)^{a-2}[1-G(t)^a]^{b-1}]}{1-\overline{\alpha}[1-G(t)^a]^b}$$

$$+\frac{\overline{\alpha}\,a\,b[\,a\,(b-1)\,g(t)^{2}\,G(t)^{2a-2}[1-G(t)^{a}\,]^{b-2}\,]}{1-\overline{\alpha}\,[1-G(t)^{a}]^{b}}+\left[\frac{\overline{\alpha}\,a\,b\,g(t)\,G(t)^{a-1}[1-G(t)^{a}\,]^{b-1}}{1-\overline{\alpha}\,[1-G(t)^{a}]^{b}}\right]^{2}$$

$$=\frac{g(t)\,g''(t)-[\,g'(t)]^{2}}{g(t)^{2}}+(a-1)\,\frac{G(t)\,g'(t)-g(t)^{2}}{G(t)}$$

$$+a\,[\,\frac{g'(t)\,G(t)^{a-1}}{1-G(t)^{a}}+\frac{(a-1)\,g(t)^{2}\,G(t)^{a-2}}{1-G(t)^{a}}+\frac{a\,g(t)^{2}\,G(t)^{2a-2}}{[1-G(t)^{a}]^{2}}\,]-\frac{\overline{\alpha}\,g'(t)}{g(t)}\,[1-G(t)^{a}\,]^{b}\,h^{MOKwG}(t)$$

$$-\overline{\alpha}\,(a-1)\,g(t)\,G(t)^{a-1}(t)\,[1-G(t)^{a}\,]^{b}\,h^{MOKwG}(t)+\overline{\alpha}\,a\,(b-1)\,g(t)\,G(t)^{a-1}(t)\,\{1-G(t)^{a}\,\}^{b-1}\,h^{MOKwG}(t)$$

$$+\{[1-G(t)^{a}\,]^{b}\,\overline{\alpha}\,h^{MOKwG}(t)\}^{2}$$

## 4.8 Stochastic orderings

In this section we study the reliability properties and stochastic ordering of the $MOKw-G$ distributions Stochastic ordering properties have applications in diverse fields such as economics, reliability, survival analysis, insurance, finance, actuarial and management sciences (Shaked and Shanthikumar, 2007).

Let $X$ and $Y$ be two random variables with cfds $F$ and $G$, respectively, survival functions $\overline{F}=1-F$ and $\overline{G}=1-G$, and corresponding pdf's $f$, $g$. Then $X$ is said to be smaller than Y in the likelihood ratio order ($X\le_{lr}Y$) if $f(t)/g(t)$ is decreasing in $t\ge 0$; stochastic order ($X\le_{st}Y$) if $\overline{F}(t)\le\overline{G}(t)$ for all $t\ge 0$; hazard rate order ($X\le_{hr}Y$) if $\overline{F}(t)/\overline{G}(t)$ is decreasing in $t\ge 0$; reversed hazard rate order ($X\le_{rhr}Y$) if $F(t)/G(t)$ is decreasing in $t\ge 0$. These four stochastic orders are related to each other, as

$$X\le_{rhr}Y\Leftarrow X\le_{lr}Y\Rightarrow X\le_{hr}Y\Rightarrow X\le_{st}Y \qquad (24)$$

**Theorem1**: Let $X\sim MOKwG(\alpha_1,a,b)$ and $Y\sim MOKwG(\alpha_2,a,b)$. If $\alpha_1<\alpha_2$, then $X\le_{lr}Y$ ($X\le_{hr}Y,\ X\le_{rhr}Y,\ X\le_{st}Y$).

Proof: $$\frac{f(t)}{g(t)}=\frac{\alpha_1\,a\,b\,g(t)\,G(t)^{a-1}[1-G(t)^{a}\,]^{b-1}/[1-\overline{\alpha}_1\,[1-G(t)^{a}\,]^{b}\,]^{2}}{\alpha_2\,a\,b\,g(t)\,G(t)^{a-1}[1-G(t)^{a}\,]^{b-1}/[1-\overline{\alpha}_2\,[1-G(t)^{a}\,]^{b}\,]^{2}}$$

$$=\frac{\alpha_1}{\alpha_2}\left[\frac{1-\overline{\alpha}_2\,[1-G(t)^{a}\,]^{b}}{1-\overline{\alpha}_1\,[1-G(t)^{a}\,]^{b}}\right]^{2}$$

Since $\alpha_1<\alpha_2$

$$\frac{d}{dt}\left[\frac{f(t)}{g(t)}\right]=2\frac{\alpha_1}{\alpha_2}\frac{1-\overline{\alpha}_2\,[1-G(t)^{a}\,]^{b}}{1-\overline{\alpha}_1\,[1-G(t)^{a}\,]^{b}}$$

$$\times\frac{\{1-\overline{\alpha}_1[1-G(t)^{a}]^{b}\}[\overline{\alpha}_2 a\,b[1-G(t)^{a}]^{b-1}G(t)^{a-1}\,g(t)]-\{1-\overline{\alpha}_2\,[1-G(t)^{a}]^{b}\}[\overline{\alpha}_1 a\,b[1-G(t)^{a}]^{b-1}G(t)^{a-1}\,g(t)]}{\{1-\overline{\alpha}_1[1-G(t)^{a}]^{b}\}^{2}}$$

$$=\frac{2\alpha_1\{1-\overline{\alpha}_2\,[1-G(t)^{a}\,]^{b}\}\,a\,b[1-G(t)^{a}\,]^{b-1}\,G(t)^{a-1}\,g(t)[\{1-\overline{\alpha}_1\,[1-G(t)^{a}\,]^{b}\,\}\overline{\alpha}_2-\{1-\overline{\alpha}_2\,[1-G(t)^{a}\,]^{b}\,\}\overline{\alpha}_1\,]}{\alpha_2\{1-\overline{\alpha}_1\,[1-G(t)^{a}\,]^{b}\,\}^{3}}$$

$$= \frac{2\alpha_1\{1-\bar{\alpha}_2[1-G(t)^a]^b\}ab[1-G(t)^a]^{b-1}G(t)^{a-1}g(t)(\bar{\alpha}_2-\bar{\alpha}_1)}{\alpha_2\{1-\bar{\alpha}_1[1-G(t)^a]^b\}^3}$$

$$= \frac{2\alpha_1(\alpha_1-\alpha_2)\{1-\bar{\alpha}_2[1-G(t)^a]^b\}ab[1-G(t)^a]^{b-1}G(t)^{a-1}g(t)}{\alpha_2\{1-\bar{\alpha}_1[1-G(t)^a]^b\}^3}.$$ Which is always less than 0.

Hence, $f(t)/g(t)$ is decreasing in *t*. That is $X \leq_{lr} Y$.The remaining statements follow from the implications (24).

## 5. Estimation

### 5.1 Maximum likelihood estimation for $MOKw-G$

The model parameters of the $MOKw-G$ distribution can be estimated by maximum likelihood. Let $t=(t_1,t_2,...t_n)^T$ be a random sample of size $n$ from $MOKw-G$ with parameter vector $\boldsymbol{\theta}=(\alpha,a,b,\boldsymbol{\beta}^T)^T$, where $\boldsymbol{\beta}=(\beta_1,\beta_2,...\beta_q)^T$ corresponds to the parameter vector of the baseline distribution *G*. Then the log-likelihood function for $\boldsymbol{\theta}$ is given by

$$\ell=\ell(\boldsymbol{\theta})$$

$$=n\log\alpha+\sum_{i=1}^{n}\log[abg(t_i,\boldsymbol{\beta})G(t_i,\boldsymbol{\beta})^{a-1}[1-G(t_i,\boldsymbol{\beta})^a]^{b-1}]-2\sum_{i=1}^{n}\log[1-\bar{\alpha}[1-G(t_i,\boldsymbol{\beta})^a]^b]$$

$$=n\log\alpha+n\log(ab)+\sum_{i=0}^{n}\log[g(t_i,\boldsymbol{\beta})]+(a-1)\sum_{i=0}^{n}\log[G(t_i,\boldsymbol{\beta})]+(b-1)\sum_{i=0}^{n}\log[1-G(t_i,\boldsymbol{\beta})^a]$$

$$-2\sum_{i=1}^{n}\log[1-\bar{\alpha}[1-G(t_i,\boldsymbol{\beta})^a]^b] \quad (25)$$

This log-likelihood function can not be solved analytically because of its complex form but it can be maximized numerically by employing global optimization methods available with softwares like R, Mathematica.

By taking the partial derivatives of the log-likelihood function with respect to $\alpha,a,b$ and $\boldsymbol{\beta}$ components of the score vector $U_{\boldsymbol{\theta}}=(U_\alpha,U_a,U_b,U_{\beta^T})^T$ can be obtained as follows:

$$U_\alpha=\frac{\partial\ell}{\partial\alpha}=\frac{n}{\alpha}-2\sum_{i=0}^{n}\frac{[1-G(t_i,\boldsymbol{\beta})^a]^b}{1-\bar{\alpha}[1-G(t_i,\boldsymbol{\beta})^a]^b}$$

$$U_a=\frac{\partial\ell}{\partial a}=\frac{n}{a}+\sum_{i=0}^{n}\log[G(t_i,\boldsymbol{\beta})]+(1-b)\sum_{i=0}^{n}\frac{G(t_i,\boldsymbol{\beta})^a\log[G(t_i,\boldsymbol{\beta})]}{1-G(t_i,\boldsymbol{\beta})^a}$$

$$-2\sum_{i=0}^{n}\frac{b\bar{\alpha}[1-G(t_i,\boldsymbol{\beta})^a]^{b-1}G(t_i,\boldsymbol{\beta})^a\log[G(t_i,\boldsymbol{\beta})]}{1-\bar{\alpha}[1-G(t_i,\boldsymbol{\beta})^a]^b}$$

$$U_b=\frac{\partial\ell}{\partial b}=\frac{n}{b}+\sum_{i=0}^{n}\log[1-G(t_i,\boldsymbol{\beta})^a]+2\sum_{i=0}^{n}\frac{\bar{\alpha}[1-G(t_i,\boldsymbol{\beta})^a]^b\log[1-G(t_i,\boldsymbol{\beta})^a]}{1-\bar{\alpha}[1-G(t_i,\boldsymbol{\beta})^a]^b}$$

$$U_{\beta}=\frac{\partial \ell}{\partial \boldsymbol{\beta}}=\sum_{i=0}^{n}\frac{g^{(\boldsymbol{\beta})}(t_i,\boldsymbol{\beta})}{g(t_i,\boldsymbol{\beta})}+(a-1)\sum_{i=0}^{n}\frac{G^{(\boldsymbol{\beta})}(t_i,\boldsymbol{\beta})}{G(t_i,\boldsymbol{\beta})}+(1-b)\sum_{i=0}^{n}\frac{aG(t_i,\boldsymbol{\beta})^{a-1}G^{(\boldsymbol{\beta})}(t_i,\boldsymbol{\beta})}{1-G(t_i,\boldsymbol{\beta})^{a}}$$

$$-2\sum_{i=0}^{n}\frac{b\overline{\alpha}[1-G(t_i,\boldsymbol{\beta})^{a}]^{b-1}aG(t_i,\boldsymbol{\beta})^{a-1}G^{(\boldsymbol{\beta})}(t_i,\boldsymbol{\beta})}{1-\overline{\alpha}\,[1-G(t_i,\boldsymbol{\beta})^{a}]^{b}}$$

Solving the equations $U_{\theta}=(U_{\alpha},U_{a},U_{b},U_{\beta^T})^T=0$ simultaneously gives the maximum likelihood estimate (MLE) $\hat{\theta}=(\hat{\alpha},\hat{a},\hat{b},\hat{\boldsymbol{\beta}}^T)^T$ of $\theta=(\alpha,a,b,\boldsymbol{\beta}^T)^T$.

**5.2 Asymptotic standard error and confidence interval for the mles**

The asymptotic variance-covariance matrix of the MLEs of parameters can obtained by inverting the Fisher information matrix $\mathrm{I}(\boldsymbol{\theta})$ which can be derived using the second partial derivatives of the log-likelihood function with respect to each parameter. The $i\,j^{th}$ elements of $\mathrm{I}_n(\boldsymbol{\theta})$ are given by

$$\mathrm{I}_{ij}=-E\left(\frac{\partial^2 l(\boldsymbol{\theta})}{\partial\theta_i\,\partial\theta_j}\right),\qquad i,\ j=1,2,\cdots,3+q$$

The exact evaluation of the above expectations may be cumbersome. In practice one can estimate $\mathrm{I}_n(\boldsymbol{\theta})$ by the observed Fisher's information matrix $\hat{\mathrm{I}}_n(\hat{\boldsymbol{\theta}})$ is defined as:

$$\hat{\mathrm{I}}_{ij}\approx\left(-\frac{\partial^2 l(\boldsymbol{\theta})}{\partial\theta_i\,\partial\theta_j}\right)_{\boldsymbol{\theta}=\hat{\boldsymbol{\theta}}},\qquad i,\ j=1,2,\cdots,3+q$$

Using the general theory of MLEs under some regularity conditions on the parameters as $n\to\infty$ the asymptotic distribution of $\sqrt{n}\,(\hat{\boldsymbol{\theta}}-\boldsymbol{\theta})$ is $N_k(0,V_n)$ where $V_n=(v_{jj})=\mathrm{I}_n^{-1}(\boldsymbol{\theta})$. The asymptotic behaviour remains valid if $V_n$ is replaced by $\hat{V}_n=\hat{\mathrm{I}}^{-1}(\hat{\boldsymbol{\theta}})$. This result can be used to provide large sample standard errors and also construct confidence intervals for the model parameters. Thus an approximate standard error and $(1-\gamma/2)100\%$ confidence interval for the mle of $j^{\text{th}}$ parameter $\theta_j$ are respectively given by $\sqrt{\hat{v}_{jj}}$ and $\hat{\theta}_j\pm Z_{\gamma/2}\sqrt{\hat{v}_{jj}}$ , where $Z_{\gamma/2}$ is the $\gamma/2$ point of standard normal distribution.

As an illustration on the MLE method its large sample standard errors, confidence interval in the case of $MOKw-E(\alpha,a,b,\lambda)$ is discussed in an appendix.

**5.3 Estimation by method of moments**

Here an alternative method to estimation of the parameters is discussed. Since the moments are not in closed form, the estimation by the usual method of moments may not be tractable. Therefore following (Barreto-Souzai *et al.*, 2013) we get

$$E[1-\bar{\alpha}\{1-\bar{\alpha}[1-G(t)^a]^b\}]=\int_{-\infty}^{\infty}1-\bar{\alpha}[1-G(t)^a]^b\frac{\alpha\,a\,b\,g(t)G(t)^{a-1}[1-G(t)^a]^{b-1}}{\left[1-\bar{\alpha}[1-G(t)^a]^b\right]^2}\,dt$$

$$=\int_{-\infty}^{\infty}\frac{\alpha\,a\,b\,g(t)G(t)^{a-1}[1-G(t)^a]^{b-1}}{1-\bar{\alpha}\,[1-G(t)^a]^b}\,dt$$

Let, $u=1-\bar{\alpha}\,[1-G(t)^a]^b$ then $du=\bar{\alpha}\,a\,b\,g(t)G(t)^{a-1}[1-G(t)^a]^{b-1}\,dt$

$$\int_{\alpha}^{1}\frac{\alpha}{\bar{\alpha}}\frac{1}{u}\,du=\frac{\alpha}{\bar{\alpha}}\left[\log u\right]_{\alpha}^{1}=-\frac{\alpha}{\bar{\alpha}}\log\alpha=-\frac{\alpha}{(1-\alpha)}\log\alpha$$

$$E\,[[1-\bar{\alpha}\{1-\bar{\alpha}[1-G(t)^a]^b\}]^v]=\int_{-\infty}^{\infty}\left\{1-\bar{\alpha}[1-G(t)^a]^b\right\}^v\frac{\alpha\,a\,b\,g(t)G(t)^{a-1}[1-G(t)^a]^{b-1}}{[1-\bar{\alpha}[1-G(t)^a]^b]^2}\,dt$$

$$=\int_{-\infty}^{\infty}\{1-\bar{\alpha}\,[1-G(t)^a]^b\}^{v-2}\alpha\,a\,b\,g(t)G(t)^{a-1}[1-G(t)^a]^{b-1}\,dt$$

$$\int_{\alpha}^{1}\frac{\alpha}{\bar{\alpha}}u^{v-2}\,du=\frac{\alpha}{\bar{\alpha}}\left[\frac{u^{v-1}}{v-1}\right]_{\alpha}^{1}=\frac{\alpha\,(1-\alpha^{v-1})}{\bar{\alpha}\,(v-1)}$$

$$E[[1-\bar{\alpha}\{1-\bar{\alpha}[1-G(t)^a]^b\}]^v]=\begin{cases}-\alpha\log(\alpha)/1-\alpha, & v=1\\ \alpha\,(1-\alpha^{v-1})/\bar{\alpha}\,(v-1), & v=(2,3,...)\end{cases}\tag{28}$$

One can use (28) to give a new method of estimation i.e. for a random sample $t_1,t_2,...t_n$ from a population with survival function (8) the model parameters can be estimated from the equations

$$\frac{1}{n}\sum_{i=1}^{n}[[1-\bar{\alpha}\left\{1-\bar{\alpha}[1-G(t_i)^a]^b\right\}]^v]=\begin{cases}-\alpha\log(\alpha)/1-\alpha, & v=1\\ \alpha\,(1-\alpha^{v-1})/\bar{\alpha}\,(v-1), & v=2,3,...q+1\end{cases}$$

## 6. Real life applications

In this subsection, we consider fitting of two real data sets to show that the distributions from the proposed $MOKw-G$ distribution can be a better model than the corresponding distributions from $KwMO-G$ (Alizadeh *et al.*, 2015) by considering the Frechet and the Exponential distribution as our $G$. Here the parameters are estimated by numerical maximization of log-likelihood function and their asymptotic standard errors and 95% confidence intervals are computed using large sample approach.

In order to compare the distributions, we have considered known criteria like AIC (Akaike Information Criterion), BIC (Bayesian Information Criterion), CAIC (Consistent Akaike Information Criterion) and HQIC (Hannan-Quinn Information Criterion)

It may be noted that $AIC=2k-2l$; $BIC=k\log(n)-2l$; $CAIC=AIC+(2k(k+1))/(n-k-1)$; and $HQIC=2k\log[\log(n)]-2l$, where $k$ is the number of parameters in the statistical model, $n$ the sample size and $l$ is the maximized value of the log-likelihood function under the considered model. We have used the maximum likelihood method for estimating the model parameters.

**Example I:**

Here we work with the following data about 346 nicotine measurements made from several brands of cigarettes in 1998. The data have been collected by the Federal Trade Commission which is an independent agency of the US government, whose main mission is the promotion of consumer protection. [http: //www.ftc.gov/ reports/tobacco or http: // pw1.netcom.com/ rdavis2/ smoke. html.]

{1.3, 1.0, 1.2, 0.9, 1.1, 0.8, 0.5, 1.0, 0.7, 0.5, 1.7, 1.1, 0.8, 0.5, 1.2, 0.8, 1.1, 0.9, 1.2, 0.9, 0.8, 0.6, 0.3, 0.8, 0.6, 0.4, 1.1, 1.1, 0.2, 0.8, 0.5, 1.1, 0.1, 0.8, 1.7, 1.0, 0.8, 1.0, 0.8, 1.0, 0.2, 0.8, 0.4, 1.0, 0.2, 0.8, 1.4, 0.8, 0.5, 1.1, 0.9, 1.3, 0.9, 0.4, 1.4, 0.9, 0.5, 1.7, 0.9, 0.8, 0.8, 1.2, 0.9, 0.8, 0.5, 1.0, 0.6, 0.1, 0.2, 0.5, 0.1, 0.1, 0.9, 0.6, 0.9, 0.6, 1.2, 1.5, 1.1, 1.4, 1.2, 1.7, 1.4, 1.0, 0.7, 0.4, 0.9, 0.7, 0.8, 0.7, 0.4, 0.9, 0.6, 0.4, 1.2, 2.0, 0.7, 0.5, 0.9, 0.5, 0.9, 0.7, 0.9, 0.7, 0.4, 1.0, 0.7, 0.9, 0.7, 0.5, 1.3, 0.9, 0.8, 1.0, 0.7, 0.7, 0.6, 0.8, 1.1, 0.9, 0.9, 0.8, 0.8, 0.7, 0.7, 0.4, 0.5, 0.4, 0.9, 0.9, 0.7, 1.0, 1.0, 0.7, 1.3, 1.0, 1.1, 1.1, 0.9, 1.1, 0.8, 1.0, 0.7, 1.6, 0.8, 0.6, 0.8, 0.6, 1.2, 0.9, 0.6, 0.8, 1.0, 0.5, 0.8, 1.0, 1.1, 0.8, 0.8, 0.5, 1.1, 0.8, 0.9, 1.1, 0.8, 1.2, 1.1, 1.2, 1.1, 1.2, 0.2, 0.5, 0.7, 0.2, 0.5, 0.6, 0.1, 0.4, 0.6, 0.2, 0.5, 1.1, 0.8, 0.6, 1.1, 0.9, 0.6, 0.3, 0.9, 0.8, 0.8, 0.6, 0.4, 1.2, 1.3, 1.0, 0.6, 1.2, 0.9, 1.2, 0.9, 0.5, 0.8, 1.0, 0.7, 0.9, 1.0, 0.1, 0.2, 0.1, 0.1, 1.1, 1.0, 1.1, 0.7, 1.1, 0.7, 1.8, 1.2, 0.9, 1.7, 1.2, 1.3, 1.2, 0.9, 0.7, 0.7, 1.2, 1.0, 0.9, 1.6, 0.8, 0.8, 1.1, 1.1, 0.8, 0.6, 1.0, 0.8, 1.1, 0.8, 0.5, 1.5, 1.1, 0.8, 0.6, 1.1, 0.8, 1.1, 0.8, 1.5, 1.1, 0.8, 0.4, 1.0, 0.8, 1.4, 0.9, 0.9, 1.0, 0.9, 1.3, 0.8, 1.0, 0.5, 1.0, 0.7, 0.5, 1.4, 1.2, 0.9, 1.1, 0.9, 1.1, 1.0, 0.9, 1.2, 0.9, 1.2, 0.9, 0.5, 0.9, 0.7, 0.3, 1.0, 0.6, 1.0, 0.9, 1.0, 1.1, 0.8, 0.5, 1.1, 0.8, 1.2, 0.8, 0.5, 1.5, 1.5, 1.0, 0.8, 1.0, 0.5, 1.7, 0.3, 0.6, 0.6, 0.4, 0.5, 0.5, 0.7, 0.4, 0.5, 0.8, 0.5, 1.3, 0.9, 1.3, 0.9, 0.5, 1.2, 0.9, 1.1, 0.9, 0.5, 0.7, 0.5, 1.1, 1.1, 0.5, 0.8, 0.6, 1.2, 0.8, 0.4, 1.3, 0.8, 0.5, 1.2, 0.7, 0.5, 0.9, 1.3, 0.8, 1.2, 0.9}

In Table 1, the MLEs and standard errors (SEs) and 95% confidence intervals (in parentheses) of the parameters from all the fitted distributions along with AIC, BIC, CAIC and HQIC values are presented. According to the lowest values of the AIC, BIC, CAIC and HQIC, the $MOKw-Fr$ ($MOKw-E$) distribution could be chosen as the better model than $KwMO-Fr$ ($KwMO-E$) distribution respectively.

More information is provided for a visual comparison is presented in the form of a histogram of the data with the fitted densities in Figures 3. The fitted cdf's are also displayed in Figures 4. These plots indicate that the proposed distributions provide a good fit to these data.

**Table 1:** MLEs, standard errors and 95% confidence intervals (in parentheses) and the AIC, BIC, CAIC and HQIC values for the nicotine measurements data.

| Parameters | $KwMO-Fr$ | $MOKw-Fr$ | $KwMO-E$ | $MOKw-E$ |
|---|---|---|---|---|
| $\hat{a}$ | 10.462<br>(2.669)<br>(5.23, 15.69) | 1.176<br>(0.497)<br>(0.20, 2.15) | 1.434<br>(0.236)<br>(0.97, 1.80) | 1.511<br>(0.515)<br>(0.50, 2.52) |
| $\hat{b}$ | 75.085<br>(19.512)<br>(36.84,113.33) | 126.659<br>(100.668)<br>(-70.65, 323.97) | 1.959<br>(0.779)<br>(0.43,3.49) | 11.063<br>(35.759)<br>(-59.02, 81.15) |
| $\hat{\lambda}$ | 0.251<br>(0.014)<br>(0.22, 0.28) | 0.303<br>(0.063)<br>(0.17, 0.43) | 3.557<br>(0.677)<br>(2.23, 4.86) | 0.586<br>(1.572)<br>(-2.49, 3.67) |
| $\hat{\delta}$ | 16.270<br>(4.834)<br>(6.79, 25.75) | 29.261<br>(23.497)<br>(-16.79,75.32) | --- | --- |
| $\hat{\alpha}$ | 0.075<br>(0.013)<br>(0.04, 0.10) | 60.442<br>(85.106)<br>(-107.37, 227.25) | 21.214<br>(10.377)<br>(0.88, 41.56) | 20.96<br>(19.978)<br>(-18.19, 60.12) |
| log-likelihood | -155.49 | **-110.27** | -107.89 | **-106.61** |
| AIC | 320.98 | **230.54** | 223.78 | **221.22** |
| BIC | 340.21 | **249.77** | 239.17 | **236.61** |
| CAIC | 321.16 | **230.72** | 223.89 | **221.34** |
| HQIC | 328.64 | **238.19** | 229.91 | **227.35** |

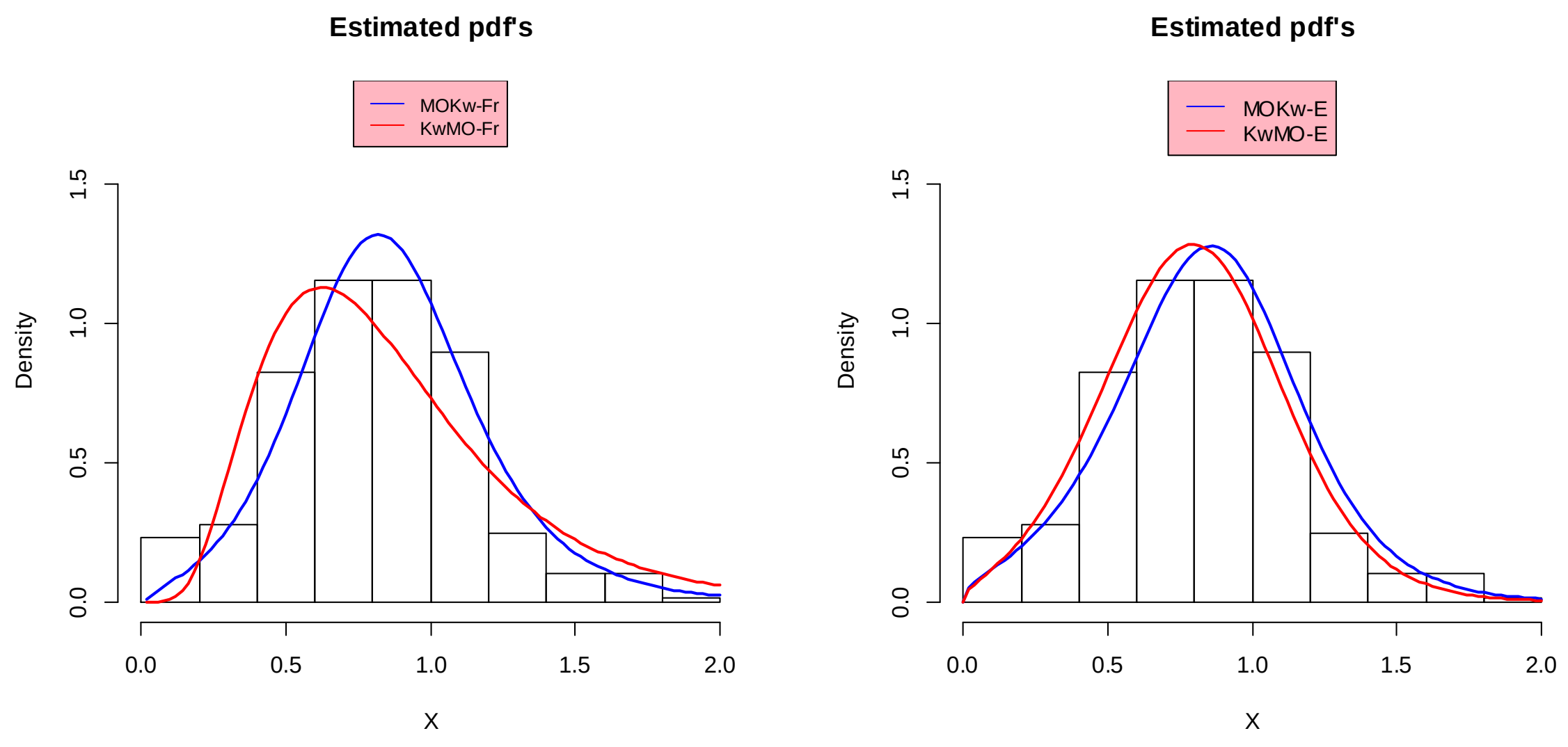

**(a)** **(b)**

**Fig: 3** Plots of the observed histogram and estimated pdf's for the (a) $MOKw-Fr$, $KwMO-Fr$ and (b) $MOKw-E$, $KwMO-E$ distributions for example I.

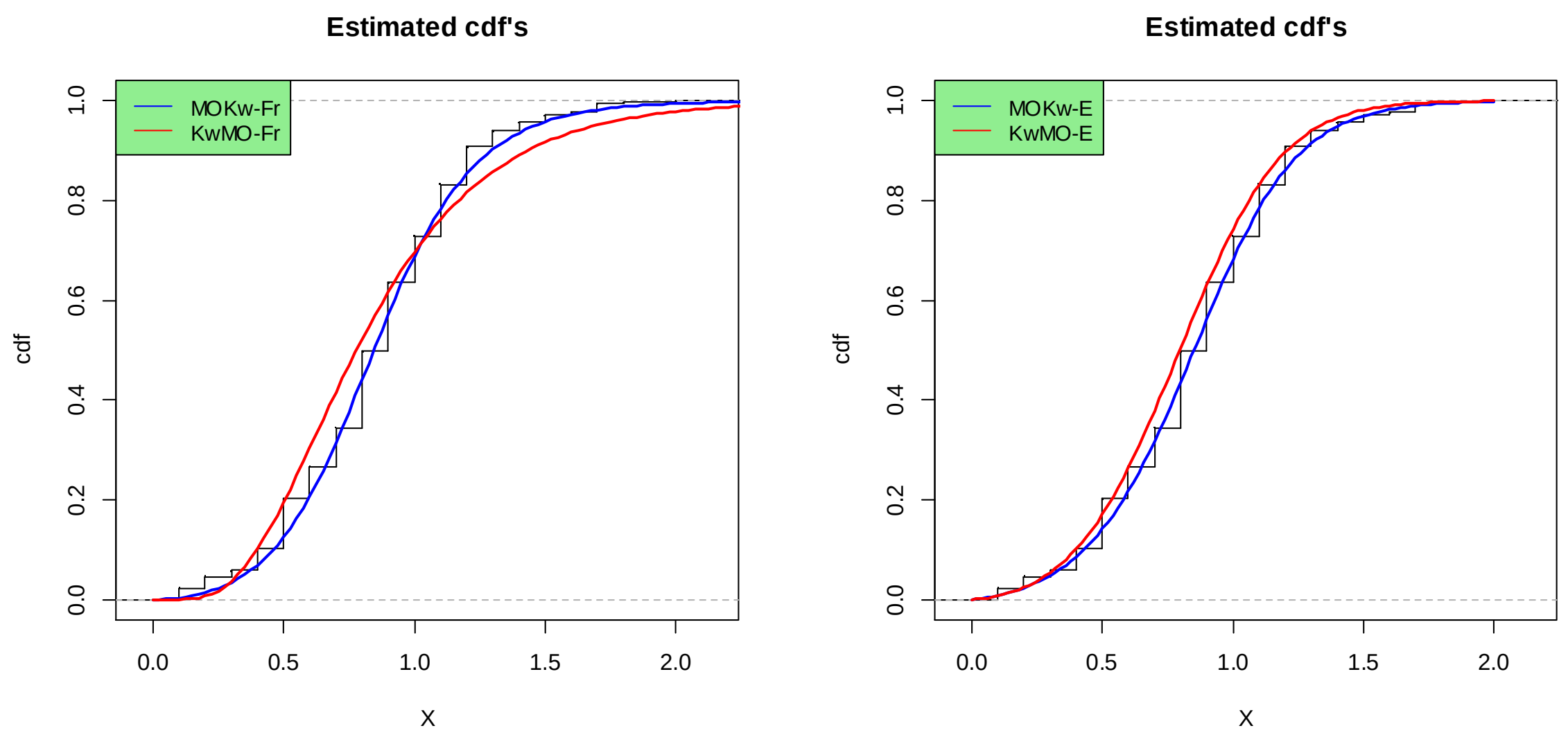

**(a)** **(b)**

**Fig: 4** Plots of the observed ogive and estimated cdf's for the (a) $MOKw-Fr$ , $KwMO-Fr$ and (b) $MOKw-E$ , $KwMO-E$ distributions for example I.

**Example II:**

This data set consists of 100 observations of breaking stress of carbon fibres (in Gba) given by Nichols and Padgett (2006).

{3.70, 2.74, 2.73, 2.50, 3.60, 3.11, 3.27, 2.87, 1.47, 3.11, 4.42, 2.40, 3.15, 2.67,3.31, 2.81, 0.98, 5.56, 5.08, 0.39, 1.57, 3.19, 4.90, 2.93, 2.85, 2.77, 2.76, 1.73, 2.48, 3.68, 1.08, 3.22, 3.75, 3.22, 2.56, 2.17, 4.91, 1.59, 1.18, 2.48, 2.03, 1.69, 2.43, 3.39, 3.56, 2.83, 3.68, 2.00, 3.51, 0.85, 1.61, 3.28, 2.95, 2.81, 3.15, 1.92, 1.84, 1.22, 2.17, 1.61, 2.12, 3.09, 2.97, 4.20, 2.35, 1.41, 1.59, 1.12, 1.69, 2.79, 1.89, 1.87, 3.39, 3.33, 2.55, 3.68, 3.19, 1.71, 1.25, 4.70, 2.88, 2.96, 2.55, 2.59, 2.97, 1.57, 2.17, 4.38, 2.03, 2.82, 2.53, 3.31, 2.38, 1.36, 0.81, 1.17, 1.84, 1.80, 2.05, 3.65}.

In Table 2, the MLEs and standard errors (SEs) and 95% confidence intervals (in parentheses) of the parameters from all the fitted distributions along with AIC, BIC, CAIC and HQIC values are presented. According to the lowest values of the AIC, BIC, CAIC and HQIC, the $MOKw-Fr$ ($MOKw-E$) distribution could be chosen as the better model than $KwMO-Fr$ ($KwMO-E$) distribution respectively.

**Table 2:** MLEs, standard errors and 95% confidence intervals (in parentheses) and the AIC, BIC, CAIC and HQIC values for the breaking stress of carbon fibres data

| Parameters | $KwMO-Fr$ | $MOKw-Fr$ | $KwMO-E$ | $MOKw-E$ |
|---|---|---|---|---|
| $\hat{a}$ | 9.243<br>(2.958)<br>(3.45, 15.04) | 6.777<br>(4.936)<br>(-2.89, 16.45) | 2.647<br>(2.066)<br>(-1.40, 6.69) | 3.226<br>(1.409)<br>(0.46, 5.99) |
| $\hat{b}$ | 19.203<br>(14.763)<br>(-9.73, 48.14) | 38.279<br>(67.729)<br>(-94.47, 171.03) | 4.571<br>(14.461)<br>(-23.77, 32.91) | 9.065<br>(42.865)<br>(-74.95, 93.08) |
| $\hat{\lambda}$ | 0.953<br>(0.172)<br>(0.62, 1.29) | 0.522<br>(0.259)<br>(0.01, 1.03) | 0.591<br>(1.103)<br>(-1.57, 0.75) | 0.295<br>(0.834)<br>(-1.34, 0.93) |
| $\hat{\delta}$ | 0.051<br>(0.026)<br>(0.00004, 0.10) | 0.416<br>(0.500)<br>(-0.56, 1.39) | --- | --- |
| $\hat{\alpha}$ | 17.222<br>(7.881)<br>(1.78, 32.67) | 16.978<br>(31.897)<br>(-45.54, 79.49) | 3.899<br>(10.290)<br>(-16.27, 24.07) | 2.566<br>(5.709)<br>(-8.62, 13.76) |
| log-likelihood | -144.83 | **-141.63** | -141.25 | **-141.09** |
| AIC | 299.66 | **293.26** | 290.50 | **290.18** |
| BIC | 312.69 | **306.29** | 300.92 | **300.60** |
| CAIC | 300.29 | **293.89** | 290.92 | **290.60** |
| HQIC | 304.94 | **298.54** | 294.73 | **294.40** |

Like for the last example here also we have presented the histogram of the observed data with the fitted densities in Figures 3 and fitted cdf's with empirical cdf in Figures 4 for visual inspection of the closeness of the fittings. These plots indicate that the distributions from the proposed family provide a closer fit to this data.

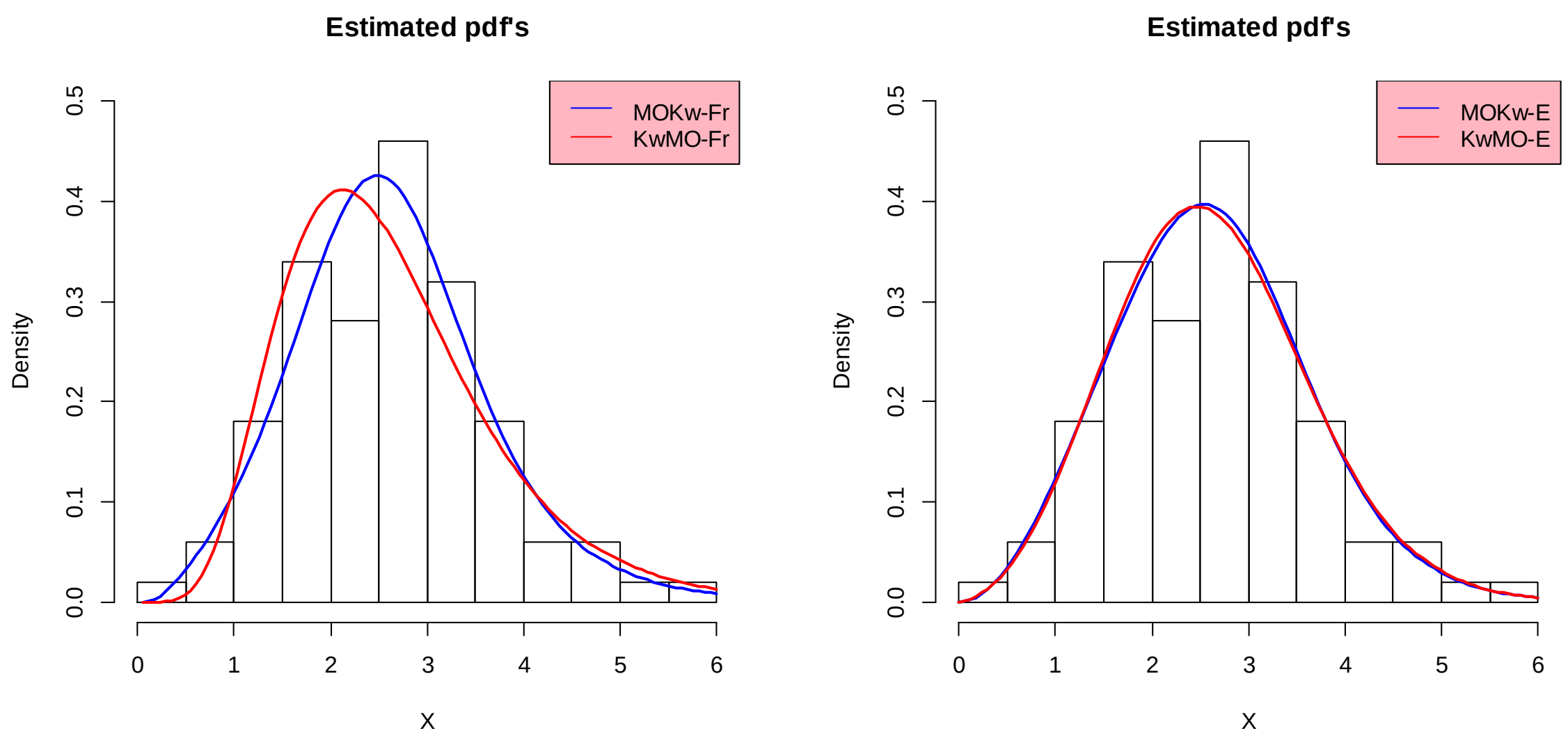

**(a)** **(b)**

**Fig: 5** Plots of the observed histogram and estimated pdf's for the (a) $MOKw-Fr$ , $KwMO-Fr$ and (b) $MOKw-E$ , $KwMO-E$ distributions for example II.

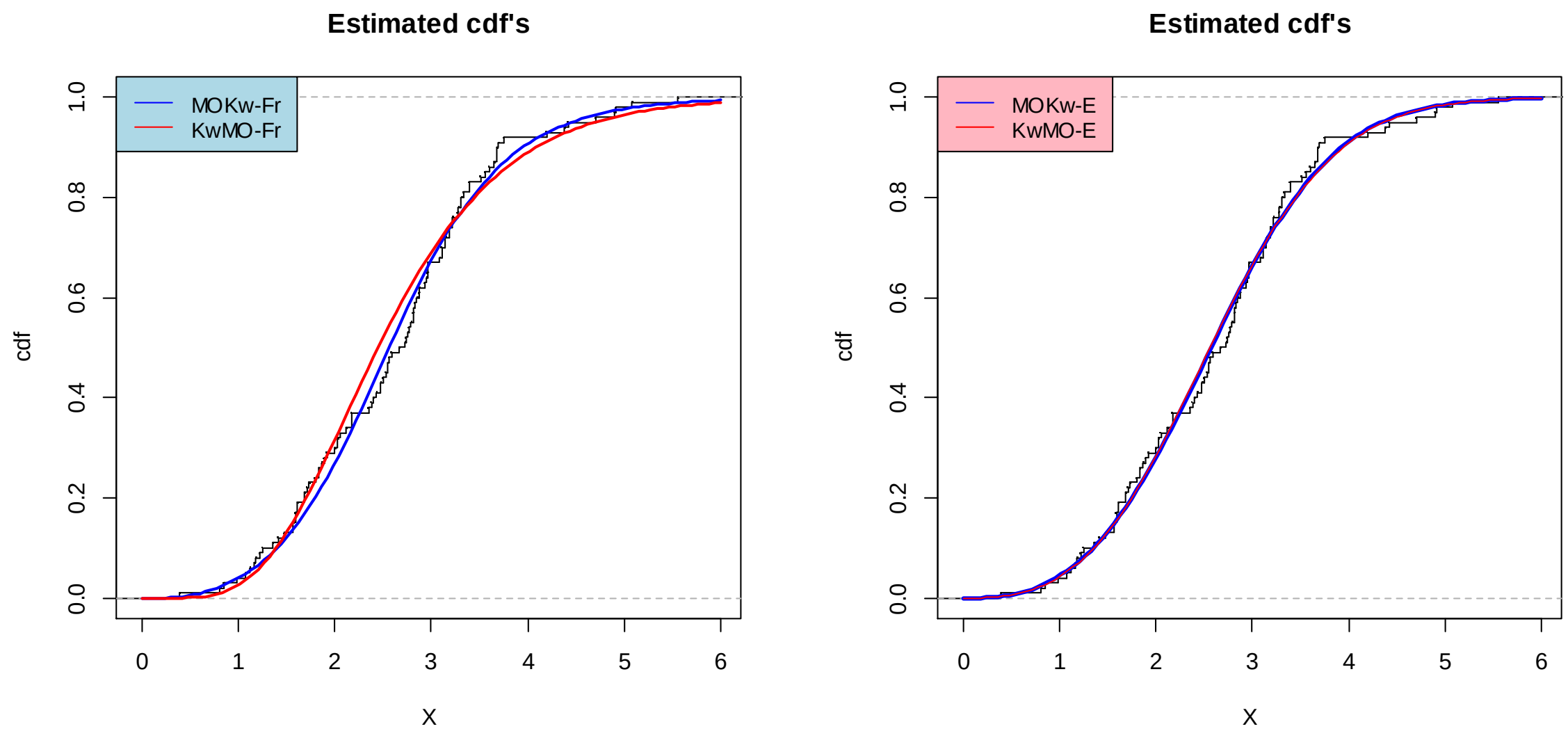

**(a)** **(b)**

**Fig: 6** Plots of the observed ogive and estimated cdf's for the (a) $MOKw-Fr$ , $KwMO-Fr$ and (b) $MOKw-E$ , $KwMO-E$ distributions for example II.

## 7. Conclusion

In this paper, the new family of Marshall-Olkin extended Kumaraswamy generalized distributions is introduced and some of its distributional and mathematical properties are investigated. The maximum likelihood and moment method for estimating the parameters are also discussed. Comparative data modeling by two important members of the proposed $MOKw-G$ family with corresponding members of $KwMO-G$ family revealed formers superiority in all examples considered here.

## Appendix: Maximum likelihood estimation for $MOKw-E$

The pdf of the $MOKw-E$ distribution is given by

$$f^{MOKwGE}(t) = \alpha\, a b \lambda e^{-\lambda t}(1-e^{-\lambda t})^{a-1}[1-(1-e^{-\lambda t})^{a}]^{b-1} \Big/ [1-\bar{\alpha}[1-(1-e^{-\lambda t})^{a}]^{b}]^{2}$$

$$\text{for } \lambda > 0, \alpha > 0,\ a > 0, b > 0,\ t > 0$$

For a random sample of size $n$ from this distribution, the log-likelihood function for the parameter vector $\boldsymbol{\theta} = (\alpha, a, b, \lambda)^T$ is given by

$$\ell(\boldsymbol{\theta}) = n\log\alpha + n\log(ab) + n\log\lambda - \lambda\sum_{i=0}^{n} t_i + (a-1)\sum_{i=0}^{n}\log(1-e^{-\lambda t_i})$$

$$+(b-1)\sum_{i=0}^{n}\log[1-(1-e^{-\lambda t_i})^{a}] \ -2\sum_{i=0}^{n}\log[1-\bar{\alpha}[1-(1-e^{-\lambda t_i})^{a}]^{b}]$$

The components of the score vector $\boldsymbol{\theta} = (\alpha, a, b, \lambda)^T$ are

$$\frac{\partial\ell(\boldsymbol{\theta})}{\partial\alpha} = \frac{n}{\alpha} - 2\sum_{i=0}^{n}\frac{[1-(1-e^{-\lambda t_i})^{a}]^{b}}{1-\bar{\alpha}[1-(1-e^{-\lambda t_i})^{a}]^{b}}$$

$$\frac{\partial\ell(\boldsymbol{\theta})}{\partial a} = \frac{n}{a} + \sum_{i=0}^{n}\log(1-e^{-\lambda t_i}) + (1-b)\sum_{i=0}^{n}\frac{(1-e^{-\lambda t_i})^{a}\log(1-e^{-\lambda t_i})}{1-(1-e^{-\lambda t_i})^{a}}$$

$$-2\sum_{i=0}^{n}\frac{\bar{\alpha}\, b[1-(1-e^{-\lambda t_i})^{a}]^{b-1}(1-e^{-\lambda t_i})^{a}\log[(1-e^{-\lambda t_i})]}{1-\bar{\alpha}[1-(1-e^{-\lambda t_i})^{a}]^{b}}$$

$$\frac{\partial\ell(\boldsymbol{\theta})}{\partial b} = \frac{n}{b} + \sum_{i=0}^{n}\log[1-(1-e^{-\lambda t_i})^{a}] + 2\sum_{i=0}^{n}\frac{\bar{\alpha}[1-(1-e^{-\lambda t_i})^{a}]^{b}\log[1-(1-e^{-\lambda t_i})^{a}]}{1-\bar{\alpha}[1-(1-e^{-\lambda t_i})^{a}]^{b}}$$

$$\frac{\partial\ell(\boldsymbol{\theta})}{\partial\lambda} = \frac{n}{\lambda} - \sum_{i=0}^{n} t_i + (a-1)\sum_{i=0}^{n}\frac{\lambda e^{-\lambda t_i}}{1-e^{-\lambda t_i}} + (1-b)\sum_{i=0}^{n}\frac{a\lambda(1-e^{-\lambda t_i})^{a-1}e^{-\lambda t_i}}{1-(1-e^{-\lambda t_i})^{a}}$$

$$-2\sum_{i=0}^{n}\frac{\bar{\alpha}\, ab\lambda[1-(1-e^{-\lambda t_i})^{a}]^{b-1}(1-e^{-\lambda t_i})^{a-1}e^{-\lambda t_i}}{1-\bar{\alpha}[1-(1-e^{-\lambda t_i})^{a}]^{b}}.$$

The asymptotic variance covariance matrix for mles of the parameters of $MOKw-E$ $(\alpha, a, b, \lambda)$ distribution is estimated by

$$\hat{\mathrm{I}}_{\mathrm{n}}{}^{-1}(\hat{\boldsymbol{\theta}}) = \begin{pmatrix} \mathrm{var}(\hat{\alpha}) & \mathrm{cov}(\hat{\alpha},\hat{a}) & \mathrm{cov}(\hat{\alpha},\hat{b}) & \mathrm{cov}(\hat{\alpha},\hat{\lambda}) \\ \mathrm{cov}(\hat{a}\hat{\alpha}) & \mathrm{var}(\hat{a}) & \mathrm{cov}(\hat{a}\hat{b}) & \mathrm{cov}(\hat{a},\hat{\lambda}) \\ \mathrm{cov}(\hat{b},\hat{\alpha}) & \mathrm{cov}(\hat{b},\hat{a}) & \mathrm{var}(\hat{b}) & \mathrm{cov}(\hat{b},\hat{\lambda}) \\ \mathrm{cov}(\hat{\lambda},\hat{\alpha}) & \mathrm{cov}(\hat{\lambda},\hat{a}) & \mathrm{cov}(\hat{\lambda},\hat{b}) & \mathrm{var}(\hat{\lambda}) \end{pmatrix}$$

Where the elements of the information matrix $\hat{\mathbf{I}}_{\mathrm{n}}(\hat{\boldsymbol{\theta}}) = \left(-\left(\frac{\partial^2 l(\boldsymbol{\theta})}{\partial\theta_i\partial\theta_j}\right)\right)_{\boldsymbol{\theta}=\hat{\boldsymbol{\theta}}}$ can be derived using the following second partial derivatives:

$$\frac{\partial^2\ell}{\partial\alpha^2} = -\frac{n}{\alpha^2} + 2\sum_{i=0}^{n}\frac{[1-(1-e^{-\lambda t_i})^a]^{2b}}{[1-\overline{\alpha}[1-(1-e^{-\lambda t_i})^a]^b]^2}$$

$$\frac{\partial^2\ell}{\partial a^2} = -\frac{n}{a^2} + (1-b)\sum_{i=0}^{n}\left(\frac{(1-e^{-\lambda t_i})^{2a}\log(1-e^{-\lambda t_i})^2}{\{1-(1-e^{-\lambda t_i})^a\}^2} + \frac{(1-e^{-\lambda t_i})^a\log(1-e^{-\lambda t_i})^2}{1-(1-e^{-\lambda t_i})^a}\right)$$

$$+2\sum_{i=0}^{n}\frac{\overline{\alpha}^2 b^2(1-e^{-\lambda t_i})^{2a}[1-(1-e^{-\lambda t_i})^a]^{2(b-1)}\log(1-e^{-\lambda t_i})^2}{\{1-\overline{\alpha}[1-(1-e^{-\lambda t_i})^a]^b\}^2}$$

$$+2\sum_{i=0}^{n}\frac{b(b-1)\overline{\alpha}\,(1-e^{-\lambda t_i})^{2a}[1-(1-e^{-\lambda t_i})^a]^{(b-2)}\log(1-e^{-\lambda t_i})^2}{1-\overline{\alpha}[1-(1-e^{-\lambda t_i})^a]^b}$$

$$-2\sum_{i=0}^{n}\frac{b\,\overline{\alpha}\,(1-e^{-\lambda t_i})^a[1-(1-e^{-\lambda t_i})^a]^{(b-1)}\log(1-e^{-\lambda t_i})^2}{1-\overline{\alpha}[1-(1-e^{-\lambda t_i})^a]^b}$$

$$\frac{\partial^2\ell}{\partial b^2} = -\frac{n}{b^2} + 2\sum_{i=0}^{n}\frac{\overline{\alpha}^2[1-(1-e^{-\lambda t_i})^a]^{2b}\log[1-(1-e^{-\lambda t_i})^a]^2}{\{1-\overline{\alpha}[1-(1-e^{-\lambda t_i})^a]^b\}^2}$$

$$+2\sum_{i=0}^{n}\frac{\overline{\alpha}\,[1-(1-e^{-\lambda t_i})^a]^b\log[1-(1-e^{-\lambda t_i})^a]^2}{1-\overline{\alpha}[1-(1-e^{-\lambda t_i})^a]^b}$$

$$\frac{\partial^2\ell}{\partial\lambda^2} = -\frac{n}{\lambda^2} + (a-1)\sum_{i=0}^{n}\left(-\frac{e^{-2\lambda t_i}t_i^2}{(1-e^{-\lambda t_i})^2} - \frac{e^{-\lambda t_i}t_i^2}{1-e^{-\lambda t_i}}\right) + (1-b)\sum_{i=0}^{n}\frac{a\,(1-e^{-\lambda t_i})^{2(a-1)}e^{-2\lambda t_i}t_i^2}{\{1-(1-e^{-\lambda t_i})^a\}^2}$$

$$+(1-b)\sum_{i=0}^{n}\frac{a(a-1)\,(1-e^{-\lambda t_i})^{a-2}e^{-2\lambda t_i}t_i^2}{1-(1-e^{-\lambda t_i})^a} - (1-b)\sum_{i=0}^{n}\frac{a\,(1-e^{-\lambda t_i})^{a-1}e^{-\lambda t_i}t_i^2}{1-(1-e^{-\lambda t_i})^a}$$

$$+2\sum_{i=0}^{n}\frac{\overline{\alpha}^2 a^2 b^2\,e^{-2\lambda t_i}(1-e^{-\lambda t_i})^{2(a-1)}[1-(1-e^{-\lambda t_i})^a]^{2(b-1)}t_i^2}{\{1-\overline{\alpha}[1-(1-e^{-\lambda t_i})^a]^b\}^2}$$

$$+2\sum_{i=0}^{n}\frac{a^2\,\overline{\alpha}\,b(b-1)e^{-2\lambda t_i}(1-e^{-\lambda t_i})^{2(a-1)}[1-(1-e^{-\lambda t_i})^a]^{b-2}t_i^2}{1-\overline{\alpha}[1-(1-e^{-\lambda t_i})^a]^b}$$

$$-2\sum_{i=0}^{n}\frac{a(a-1)\overline{\alpha}\,b\,e^{-2\lambda t_i}(1-e^{-\lambda t_i})^{a-2}[1-(1-e^{-\lambda t_i})^a]^{b-1}t_i^2}{1-\overline{\alpha}[1-(1-e^{-\lambda t_i})^a]^b}$$

$$+2\sum_{i=0}^{n}\frac{a\,b\bar{\alpha}\;e^{-\lambda t_i}(1-e^{-\lambda t_i})^{a-1}[1-(1-e^{-\lambda t_i})^{a}]^{b-1}t_i^{\,2}}{1-\bar{\alpha}[1-(1-e^{-\lambda t_i})^{a}]^{b}}$$

$$\frac{\partial^2\ell}{\partial\alpha\,\partial a}=2\sum_{i=0}^{n}\frac{b\bar{\alpha}(1-e^{-\lambda t_i})^{a}[1-(1-e^{-\lambda t_i})^{a}]^{2b-1}\log(1-e^{-\lambda t_i})}{[1-\bar{\alpha}[1-(1-e^{-\lambda t_i})^{a}]^{b}]^{2}}$$

$$+2\sum_{i=0}^{n}\frac{b\,(1-e^{-\lambda t_i})^{a}[1-(1-e^{-\lambda t_i})^{a}]^{b-1}\log(1-e^{-\lambda t_i})}{1-\bar{\alpha}[1-(1-e^{-\lambda t_i})^{a}]^{b}}$$

$$\frac{\partial^2\ell}{\partial\alpha\,\partial b}=-2\sum_{i=0}^{n}\frac{\bar{\alpha}[1-(1-e^{-\lambda t_i})^{a}]^{2b}\log[1-(1-e^{-\lambda t_i})^{a}]}{\{1-\bar{\alpha}[1-(1-e^{-\lambda t_i})^{a}]^{b}\}^{2}}-2\sum_{i=0}^{n}\frac{[1-(1-e^{-\lambda t_i})^{a}]^{b}\log[1-(1-e^{-\lambda t_i})^{a}]}{1-\bar{\alpha}[1-(1-e^{-\lambda t_i})^{a}]^{b}}$$

$$\frac{\partial^2\ell}{\partial\alpha\,\partial\lambda}=2\sum_{i=0}^{n}\frac{ab\bar{\alpha}\,e^{-\lambda t_i}\,[1-(1-e^{-\lambda t_i})^{a}]^{2b-1}(1-e^{-\lambda t_i})^{a-1}t_i}{[1-\bar{\alpha}[1-(1-e^{-\lambda t_i})^{a}]^{b}]^{2}}$$

$$+2\sum_{i=0}^{n}\frac{ab\,e^{-\lambda t_i}\,[1-(1-e^{-\lambda t_i})^{a}]^{b-1}(1-e^{-\lambda t_i})^{a-1}t_i}{1-\bar{\alpha}[1-(1-e^{-\lambda t_i})^{a}]^{b}}$$

$$\frac{\partial^2\ell}{\partial a\,\partial b}=-\sum_{i=0}^{n}\frac{(1-e^{-\lambda t_i})^{a}\log(1-e^{-\lambda t_i})}{1-(1-e^{-\lambda t_i})^{a}}-2\sum_{i=0}^{n}\frac{\bar{\alpha}(1-e^{-\lambda t_i})^{a}[1-(1-e^{-\lambda t_i})^{a}]^{b-1}\log(1-e^{-\lambda t_i})}{1-\bar{\alpha}[1-(1-e^{-\lambda t_i})^{a}]^{b}}$$

$$-2\sum_{i=0}^{n}\frac{b\bar{\alpha}^{2}(1-e^{-\lambda t_i})^{a}[1-(1-e^{-\lambda t_i})^{a}]^{2b-1}\log(1-e^{-\lambda t_i})\log[1-(1-e^{-\lambda t_i})^{a}]}{\{1-\bar{\alpha}[1-(1-e^{-\lambda t_i})^{a}]^{b}\}^{2}}$$

$$-2\sum_{i=0}^{n}\frac{b\bar{\alpha}(1-e^{-\lambda t_i})^{a}[1-(1-e^{-\lambda t_i})^{a}]^{b-1}\log(1-e^{-\lambda t_i})\log[1-(1-e^{-\lambda t_i})^{a}]}{1-\bar{\alpha}[1-(1-e^{-\lambda t_i})^{a}]^{b}}$$

$$\frac{\partial^2\ell}{\partial a\,\partial\lambda}=\sum_{i=0}^{n}\frac{e^{-\lambda t_i}t_i}{1-e^{-\lambda t_i}}+(1-b)\sum_{i=0}^{n}\frac{a\,(1-e^{-\lambda t_i})^{2a-1}e^{-\lambda t_i}t_i\log(1-e^{-\lambda t_i})}{\{1-(1-e^{-\lambda t_i})^{a}\}^{2}}$$

$$+(1-b)\sum_{i=0}^{n}\frac{(1-e^{-\lambda t_i})^{a-1}e^{-\lambda t_i}t_i}{1-(1-e^{-\lambda t_i})^{a}}+(1-b)\sum_{i=0}^{n}\frac{a\,(1-e^{-\lambda t_i})^{a-1}e^{-\lambda t_i}t_i\log(1-e^{-\lambda t_i})}{1-(1-e^{-\lambda t_i})^{a}}$$

$$-2\sum_{i=0}^{n}\frac{b\bar{\alpha}\;e^{-\lambda t_i}(1-e^{-\lambda t_i})^{a-1}[1-(1-e^{-\lambda t_i})^{a}]^{b-1}t_i}{1-\bar{\alpha}[1-(1-e^{-\lambda t_i})^{a}]^{b}}$$

$$+2\sum_{i=0}^{n}\frac{a\,b^{2}\,\bar{\alpha}^{2}\,e^{-\lambda t_i}(1-e^{-\lambda t_i})^{2a-1}[1-(1-e^{-\lambda t_i})^{a}]^{2(b-1)}t_i\log(1-e^{-\lambda t_i})}{\{1-\bar{\alpha}[1-(1-e^{-\lambda t_i})^{a}]^{b}\}^{2}}$$

$$+2\sum_{i=0}^{n}\frac{a\,b(b-1)\bar{\alpha}\,e^{-\lambda t_i}(1-e^{-\lambda t_i})^{2a-1}[1-(1-e^{-\lambda t_i})^{a}]^{b-2}t_i\log(1-e^{-\lambda t_i})}{1-\bar{\alpha}[1-(1-e^{-\lambda t_i})^{a}]^{b}}$$

$$-2\sum_{i=0}^{n}\frac{a\,b\,\bar{\alpha}\,e^{-\lambda t_i}(1-e^{-\lambda t_i})^{a-1}[1-(1-e^{-\lambda t_i})^{a}]^{b-1}t_i\log(1-e^{-\lambda t_i})}{1-\bar{\alpha}[1-(1-e^{-\lambda t_i})^{a}]^{b}}$$

$$\frac{\partial^2\ell}{\partial b\,\partial\lambda}=-\sum_{i=0}^{n}\frac{ae^{-\lambda t_i}(1-e^{-\lambda t_i})^{a-1}t_i}{1-(1-e^{-\lambda t_i})^{a}}-2\sum_{i=0}^{n}\frac{a\bar{\alpha}\;e^{-\lambda t_i}(1-e^{-\lambda t_i})^{a-1}[1-(1-e^{-\lambda t_i})^{a}]^{b-1}t_i}{1-\bar{\alpha}[1-(1-e^{-\lambda t_i})^{a}]^{b}}$$

$$-2\sum_{i=0}^{n}\frac{a\,b\,\overline{\alpha}^{\,2}\,e^{-\lambda t_i}(1-e^{-\lambda t_i})^{a-1}[1-(1-e^{-\lambda t_i})^{a}]^{2b-1}t_i\log[1-(1-e^{-\lambda t_i})^{a}]}{\{1-\overline{\alpha}[1-(1-e^{-\lambda t_i})^{a}]^{b}\}^{2}}$$

$$-2\sum_{i=0}^{n}\frac{a\,b\,\overline{\alpha}\,e^{-\lambda t_i}(1-e^{-\lambda t_i})^{a-1}[1-(1-e^{-\lambda t_i})^{a}]^{b-1}t_i\log[1-(1-e^{-\lambda t_i})^{a}]}{1-\overline{\alpha}[1-(1-e^{-\lambda t_i})^{a}]^{b}}$$

Where $\psi'(.)$ is the derivative of the digamma function.